\documentclass[onefignum,onetabnum]{siamart190516}
\usepackage{lipsum}
\usepackage{amsfonts}
\usepackage{graphicx,epstopdf}
\usepackage{ntheorem}
\usepackage{verbatim}

\usepackage{booktabs}
\usepackage{algorithm,algorithmic}
\usepackage{url}
\graphicspath{{./pics/}}

\newsiamremark{remark}{Remark}

\newtheorem{example}{Example}[section]

\newtheorem{defn}{Definition}[section]
\numberwithin{equation}{section}

\usepackage{booktabs}
\usepackage{threeparttable}

\renewcommand{\d}{\mathrm{d}}
\allowdisplaybreaks
\providecommand{\abs}[1]{\left\lvert#1\right\rvert}
\providecommand{\norm}[1]{\left\lVert#1\right\rVert}

\title{Recovery of the Order of Derivation for Fractional Diffusion Equations in an Unknown Medium\thanks{The work of B.J. is partially supported by UK EPSRC grant EP/T000864/1, and that of Y.K. by the French National Research Agency ANR
(project MultiOnde) grant ANR-17-CE40-0029.}}

\author{Bangti Jin\thanks{Department of Computer Science, University College London, Gower Street, London WC1E 6BT, UK (b.jin@ucl.ac.uk)}\and
Yavar Kian\thanks{Aix Marseille Universit\'{e}, Universit\'{e} de Toulon, CNRS, CPT, Marseille, France (yavar.kian@univ-amu.fr)}}

\date{}

\begin{document}
\maketitle

\begin{abstract}
In this work, we investigate the recovery of a parameter in a diffusion process given by  the order of
derivation in time for a class of diffusion type equations, including both classical and  time-fractional diffusion equations,
from the flux measurement observed at one point on the boundary. The mathematical model for time-fractional
diffusion equations involves a Djrbashian-Caputo fractional derivative in time. We prove a uniqueness
result in an unknown medium  (e.g., diffusion coefficients, obstacle, initial condition and source), i.e.,
the recovery of the order of derivation in a diffusion process having several pieces of
unknown information. The proof relies on the analyticity of the solution at large time, asymptotic
decay behavior, strong maximum principle of the elliptic problem and suitable application of the Hopf lemma.
Further we provide an easy-to-implement reconstruction algorithm based on a nonlinear least-squares formulation, and several
numerical experiments are presented to complement the theoretical analysis.
\end{abstract}



\section{Introduction}
Let $\widetilde{\Omega}\subset \mathbb{R}^d$ ($d\ge 2$) be an open bounded and connected subset with a
$ C^{2\lceil\frac{d}{4}\rceil+2}$ boundary $\partial \widetilde{\Omega}$ (with $\lceil\cdot\rceil$ being
the ceiling function), $\omega$ a $C^{2\lceil\frac{d}{4}\rceil+2}$ open set of $\mathbb{R}^d$
such that $\overline{\omega}\subset \widetilde{\Omega}$, and let $\Omega=\widetilde{\Omega}\setminus\overline{\omega}$.
We denote by  $\nu(x)$ the unit outward normal vector to the (outer) boundary $\partial\widetilde{\Omega}$
at a point $x\in\partial\widetilde\Omega$,
and $\partial_\nu$ the normal derivative. Next we define an elliptic operator $\mathcal{A}$  on the domain $\Omega$ by
\begin{equation}\label{A}
\mathcal{A} u(x) :=-\sum_{i,j=1}^d \partial_{x_i} \left( a_{i,j}(x) \partial_{x_j} u(x) \right)+q(x)u(x),\quad  x\in\Omega,
\end{equation}
where the potential $q \in C^{2\lceil\frac{d}{4}\rceil}(\overline{\Omega})$ is nonnegative,
and the diffusion coefficient matrix $a:=(a_{i,j})_{1 \leq i,j \leq d} \in C^{1+2\lceil\frac{d}{4}\rceil}
(\overline{\Omega};\mathbb{R}^{d\times d})$ is symmetric 
and fulfills the following ellipticity condition
\begin{equation}\label{ell}
\exists c>0,\ \sum_{i,j=1}^d a_{i,j}(x) \xi_i \xi_j \geq c |\xi|^2,\quad \forall x \in \overline{\Omega},\ \forall \xi=(\xi_1,\ldots,\xi_d) \in \mathbb{R}^d.
\end{equation}
Let $\rho \in C^{2\lceil\frac{d}{4}\rceil}(\overline{\Omega})$ obey that for some $\rho_1>\rho_0>0$,
\begin{equation}\label{eqn:rho}
 0<\rho_0 \leq\rho(x) \leq\rho_1 <+\infty\quad\mbox{in } \Omega.
\end{equation}
For $\alpha\in(0,2)$ and $0<T<+\infty$, consider the following initial boundary value problem for $u$:
\begin{equation}\label{eq1}
\begin{cases}
\rho(x)\partial_t^{\alpha}u +\mathcal{A} u =  F, & \mbox{in }\Omega\times(0,T),\\
u= g, & \mbox{on } \partial{\widetilde\Omega}\times(0,T), \\
 u=0,& \mbox{on }\partial\omega\times(0,T),\\
\begin{cases}
u=u_0 & \mbox{if }0<\alpha\leq1,\\
u=u_0,\quad \partial_t u=0 & \mbox{if }1<\alpha<2,
\end{cases} & \mbox{in }\Omega\times \{0\}.
\end{cases}
\end{equation}
In the model \eqref{eq1}, the notation $\partial_t^\alpha u$ denotes the so-called
Djrbashian-Caputo fractional derivative of order $\alpha$ with respect to
$t$, which, for $\alpha \in (0,1) \cup (1,2)$, is defined by
\cite{KilbasSrivastavaTrujillo:2006,P,Jin:2021book}
\begin{equation}\label{cap}
 \partial_t^\alpha u(x,t):=\frac{1}{\Gamma(\lceil\alpha\rceil-\alpha)}\int_0^t(t-s)^{\lceil\alpha\rceil-1-\alpha}\partial_s^{\lceil\alpha\rceil} u(x,s) \d s,\quad (x,t) \in \Omega\times(0,T),
\end{equation}
where the notation $\Gamma(z)=\int_0^\infty s^{z-1}e^{-s}\d s$, $\Re(z)>0$, denotes Euler's Gamma function,
whereas for $\alpha=1$, $\partial_t^\alpha u$ is identified with the usual first order derivative
$\partial_t u$.
Throughout we assume that there exists some $T_1\in(0,T)$ such that
\begin{align}
   F(x,t)&=0,\quad (x,t)\in\Omega\times(T_1,T), \label{eqn:source}\\
   g(x,t)& =0, \quad (x,t)\in\partial\widetilde\Omega\times(T_1,T).\label{eqn:diri-bc}
\end{align}
Note that the conditions \eqref{eqn:source}--\eqref{eqn:diri-bc} require the
source $F$ and Dirichlet input $g$ vanish for the time interval $(T_1,T)$. These conditions
are needed to ensure the analyticity of the solution $u(x,t)$ in time $t$ for any $t\in(T_1,T)$, and play an
essential role in the proof of Theorem \ref{thm:main}.

The model \eqref{eq1}, with $\alpha\neq1$, is widely employed to describe anomalous diffusion processes arising
in physics, engineering and biology. The cases $\alpha\in(0,1)$ and $\alpha\in(1,2)$ are known as
subdiffusion and diffusion wave, respectively. The former can be viewed as the macroscopic counterpart of continuous
time random walk with a waiting time distribution being heavy tailed (i.e., with a divergent mean)
{in the sense that the probability density function of the particle appearing at time $t>0$ and spatial location
$x\in\mathbb{R}^d$ satisfies a differential equation of the form \eqref{eq1}.} Subdiffusion has been
observed in diffusion in media with fractal geometry \cite{Nigmatulin:1986},
transport in column experiments \cite{HatanoHatano:1998} and subsurface flows \cite{AdamsGelhar:1992} etc, whereas the
diffusion wave case was employed in dynamic viscoelasticity, describing the propagation of mechanical diffusive waves
in viscoelastic media which exhibit a power-law creep \cite{Mainardi:1995,Mainardi:2010}.
We refer interested readers to \cite{MetzlerKlafter:2000} for physical motivations and many applications.

This paper is concerned with the following inverse problem: to determine the order $\alpha$ of
the fractional derivative $\partial_t^\alpha u(x,t)$ in the model \eqref{eq1} from a knowledge
of the flux data
\begin{equation*}
  h(t):=\partial_\nu u(x_0,t),\quad  t\in (T-\delta,T),
\end{equation*}
for an arbitrary point $x_0\in\partial\widetilde{\Omega}$ and $\delta\in(0,T-T_1)$, where $u$ is a
solution to problem \eqref{eq1}, but without assuming a full knowledge of the problem data
(e.g., $u_0$, $F$, $g$, $\omega$, $a$ and $q$) in the model \eqref{eq1}. Note that the obstacle
$\omega$ is contained in the domain $\widetilde{\Omega}$,
and the direct problem \eqref{eq1} is posed on the domain $\Omega=\widetilde{\Omega}\setminus\overline{\omega}$
with a boundary $\partial\Omega=\partial\widetilde{\Omega}\cup\partial{\omega}$.  We impose a zero
Dirichlet boundary condition on $\partial\omega$ (i.e., solid obstacle), but allow a more general
Dirichlet input $g$ on $\partial\widetilde{\Omega}$. The measurement of the flux $\partial_\nu
u(x_0,t)$ at $x_0\in \partial\widetilde{\Omega}$ for the time interval $t\in(T-\delta,T)$ is performed
on the part of the boundary $\partial\Omega$ not intersecting the obstacle $\omega$, which represents
the overposed data for order determination.

The determination of fractional order(s)  is probably one of the most
natural inverse problems for time-fractional models, as was recently highlighted
by the survey \cite{LiLiuYamamoto:2019}; see also \cite{JinRundell:2015} for a tutorial on inverse
problems for anomalous diffusion. The determination of this parameter allows one to distinguish the
type of the concerned diffusion phenomenon, i.e., a classical one (corresponding to the case $\alpha=1$) or an anomalous
one described by a subdiffusive ($\alpha\in(0,1)$) or a superdiffusive ($\alpha\in(1,2)$)  model.
For subdiffusion, this inverse problem was first studied by Hatano et al \cite{HNWY}, which provided
two reconstruction formulas, based on the asymptotics of the solution at small or large time, respectively,
and also discussed the numerical recovery for smooth observational data. A first Lipschitz stability
result was recently shown {in} \cite{LiHuangYamamoto:2020}. The work \cite{Alimov:2020} gave a uniqueness
result from the terminal measurement. See also \cite{TatarTinazepeUlusoy:2016} for numerical recovery.
Krasnoschok et al \cite{KrasnoschokPereveryev:2020} studied the recovery of the order in semilinear
subdiffusion. There are several works on the simultaneous recovery of the order with the source or other
unknowns \cite{J,JK,LiZhang:2020,Lukashchuk:2011,JinZhou:2021ip}. We also refer readers to the order recovery in more
complex models, e.g., multiple orders \cite{LiYamamoto:2015,JinKian:2021multi}, spatially-variable order
\cite{KianSoccorsiYamamoto:2018}, weight in distributed-order \cite{RundellZhang:2017,LiFujishiroLi:2020}
and time-variable order \cite{ZhengChengWang:2019}. Finally, we mention \cite{TatarUlusoy:2013,Yamamoto:2020order}
dealing with similar problems for space-time fractional models and \cite{Kian:2020single} on the simultaneous
recovery of the order of derivation with coefficients, a source term and an obstacle.

In this work we consider solutions of problem \eqref{eq1} in the following sense.
\begin{defn}
\label{d3}
A function $u\in W^{\lceil\alpha\rceil,1}(0,T;H^{-1}(\Omega))\cap L^1(0,T;H^1(\Omega))$
is said to be a solution to \eqref{eq1} if $u$ solves
$\rho(x) \partial_t^{\alpha}u+\mathcal A u=F$ in the sense of $L^1(0,T;H^{-1}(\Omega))$ and satisfies
\begin{equation*}\begin{cases}
u= g, & \mbox{on } \partial{\widetilde\Omega}\times(0,T), \\
 u=0,& \mbox{on }\partial\omega\times(0,T),\\
\begin{cases}
u=u_0 & \mbox{if }0<\alpha\leq1,\\
u=u_0,\quad \partial_t u=0 & \mbox{if }1<\alpha<2,
\end{cases} & \mbox{in }\Omega\times \{0\}.
\end{cases}
\end{equation*}
\end{defn}

That is, the governing equation holds in the sense of distribution in $\Omega\times (0,T)$,
and the initial and boundary conditions are in the sense of traces of functions
$u\in W^{\lceil\alpha\rceil,1}(0,T;H^{-1}(\Omega))\cap L^1(0,T;H^1(\Omega))$.
In view of \cite[Theorems 2.5 and 2.7, Prop. 2.6 and 2.9]{KY3}, the regularity
conditions $u_0\in L^2(\Omega)$, $g\in \cap_{k=0}^{\lceil\alpha\rceil} W^{\lceil\alpha\rceil-k,1}(0,T; H^{\frac{1}{2}+\frac{k}{\lceil\alpha\rceil}}(\partial\widetilde{\Omega}))$ and
$F\in W^{\lceil\alpha\rceil,1}(0,T;L^2(\Omega))$ along with {the following conditions at $t=0$}:
$\partial_t^kg(\cdot,0)\equiv0$, $\partial_t^kF(\cdot,0)\equiv0$ for $k=0,\ldots,\lceil\alpha\rceil-1,$
imply the existence of a solution of \eqref{eq1} in the sense of Definition \ref{d3}, whereas the
uniqueness follows from \cite[Theorems 2.1--2.4]{SY}, and \cite[Theorem 4.2, Chapter 4]{LM2} for $\alpha=1$.
{We restrict the discussions to solutions in the sense of Definition \ref{d3}
for the ease of exposition. Note that for $\alpha\in(0,1)$, the unique existence of a weak solution
does not require the conditions at $t=0$ and holds under weaker regularity conditions on the problem
data (see, e.g., \cite{SY,Zacher:2009} and \cite[Chapter 4]{KubicaYamamoto:2020}). Further, the concept of weak solutions may be stated
by means of Laplace transform \cite[Theorem 2.3]{KY3}, for which the conditions at $t=0$
are not needed and the regularity on the problem data can be relaxed. Theorem \ref{thm:main}
remains valid under these weaker conditions.}


 Now we give the main (regularity) assumptions on the
problem data. For an admissible tuple, there exists a unique  solution to problem \eqref{eq1}.
\begin{defn}
A tuple $(\alpha,\omega,a,q,\rho,u_0,F,g)$ is said to be admissible if the following conditions are fulfilled for $p>\frac d2$.
\begin{itemize}
\item[{\rm(i)}] $\alpha\in(0,2)$, $\omega\subset \mathbb{R}^d$ is a $C^{2\lceil\frac{d}{4}\rceil+2}$ open set
such that $\overline{\omega}\subset \widetilde{\Omega}$, $\Omega=\widetilde\Omega\setminus\overline{\omega}$.
\item[{\rm(ii)}] $a:=(a_{i,j})_{1 \leq i,j \leq d} \in C^{1+2\lceil\frac{d}{4}\rceil}(\overline{\Omega};
\mathbb{R}^{d\times d})$ satisfies the ellipticity condition \eqref{ell}, $q \in C^{2\lceil\frac{d}{4}\rceil}
(\overline{\Omega})$ is nonnegative, $\rho \in C^{2\lceil\frac{d}{4}\rceil}(\overline{\Omega})$
obeys the condition \eqref{eqn:rho}.
\item[{\rm(iii)}] $u_0\in L^{2p}(\Omega)$, $g\in 
L^1(0,T; W^{2-\frac{1}{p},p}(\partial\widetilde{\Omega}))\cap_{k=0}^{\lceil\alpha\rceil} W^{\lceil\alpha\rceil-k,1}(0,T; H^{\frac{1}{2}+\frac{k}{\lceil\alpha\rceil}}(\partial\widetilde{\Omega}))$ satisfies \eqref{eqn:diri-bc},
$F\in L^1(0,T;L^{p}(\Omega)) \cap W^{\lceil\alpha\rceil,1}(0,T;L^2(\Omega))$ satisfies \eqref{eqn:source}, and
$\partial_t^kg(\cdot,0)\equiv0$, $\partial_t^kF(\cdot,0)\equiv0$ for $k=0,\ldots,\lceil\alpha\rceil-1$.
\end{itemize}
\end{defn}

Now we can state the main theoretical result.
\begin{theorem}\label{thm:main}
Let $(\alpha_k,\omega_k,a^k,q_k,\rho_k,u_0^k,F_k,g_k)$, $k=1,2$, be two admissible tuples, $u^k$, $k=1,2$,
be the corresponding  solution of problem \eqref{eq1} on {the domain $\Omega_k=\widetilde
\Omega\setminus\overline{\omega}_k$,} and one of the following
conditions be fulfilled
\begin{itemize}
\item[$\rm(i)$] $u_0^k\not\equiv0$, $k=1,2$, is of constant sign.
\item[$\rm(ii)$]  $u_0^k\equiv0$,  {for $F_k^*=\int_0^TF_k(t)\d t$ and $g_k^*=\int_0^Tg_k(t)\d t$ we have either $g_k^*\geq0$ and $F_k^*\geq0$ {\rm(}or $g_k^*\leq0$ and $F_k^*\leq0${\rm)}, $k=1,2$.}  Moreover,  $F_k^*\not\equiv0$ or $g_k^*\not\equiv0$, $k=1,2$.
\end{itemize}
Then, for any $0\leq T_1<T$ and $\delta\in(0,T-T_1)$, we have $u^k\in C([T-\delta,T]; C^1(\overline{\Omega_k}))$, $k=1,2.$
Moreover, for any arbitrarily chosen $\delta\in(0,T-T_1)$ and $x_0\in\partial\widetilde{\Omega}{\subset\partial\Omega_1\cap \partial\Omega_2}$, the condition
\begin{equation}\label{t1a}
\partial_\nu u^1(x_0,t)=\partial_\nu u^2(x_0,t),\quad\forall t\in(T-\delta,T)
\end{equation}
implies $\alpha_1=\alpha_2$.
\end{theorem}

The result of Theorem \ref{thm:main} is independent of the choice of the problem data $\omega_k$, $a^k$,
$q_k $, $\rho_k$, $u_0^k$, $F_k$ and $g_k$, $k=1,2$, so long as they satisfy suitable mild assumptions, i.e.,
condition \eqref{ell}-\eqref{eqn:rho}, \eqref{eqn:source} / \eqref{eqn:diri-bc}, and one of the conditions
(i), (ii). Thus, Theorem \ref{thm:main} still holds even if $\omega_1\neq\omega_2$, $a^1\neq a^2$, $q_1\neq q_2$,
$\rho_1\neq\rho_2$, $u_0^1\neq u_0^2$, $F_1\neq F_2$ and $g_1\neq g_2$, i.e., corresponding to the unique
recovery of the fractional order $\alpha$ in an unknown medium, due to the possibly unknown problem data
$\{\omega,a,q,\rho,u_0,F,g\}$ in the model \eqref{eq1}. In addition, we develop an algorithm for recovering
the fractional order $\alpha$ based on a nonlinear least-squares formulation, and illustrate the
feasibility of the approach on several one- and two-dimensional
numerical tests. The numerical results show that subdiffusion and diffusion wave exhibit distinctly different
features for the numerical recovery.

To the best of our knowledge, Theorem \ref{thm:main} is the first result on the recovery of the order of
derivation for time-fractional models in an unknown medium from a point measurement. It also seems that Theorem \ref{thm:main} is the first
result of this type stated with a Neumann boundary measurement at an arbitrary point on the boundary of
the domain $\Omega$. Indeed, in all existing results that we are aware of the medium is always assumed to be known and, in most
of these results, the measurement corresponds to the Dirichlet trace of solutions at one internal point
(see e.g. \cite{HNWY,Yamamoto:2020order}). Note that the measurements at one internal point require at least some
\textit{a priori} knowledge of the medium that can be removed while considering boundary measurement. The result
of Theorem \ref{thm:main} can for instance be applied to the recovery of the order of derivation in time
in a diffusion process for which several pieces of information (e.g. density of the medium, source of diffusion,
location of an obstacle...) are unknown.

The key tools in the analysis include smoothing properties and analyticity in time of
the solution $u$ (or its extension $\tilde u$) of problem \eqref{eq1} for large time; see Propositions \ref{p1}
and \ref{p2}. These properties are derived from a new solution representation, asymptotics
of Mittag-Leffler functions and properties of elliptic regularization. The adopted proof techniques allow us
to state  the main result for a large class of source terms $F$ and initial conditions $u_0$ by only assuming
$F\in L^1(0,T;L^{p}(\Omega)) \cap W^{\lceil\alpha\rceil,1}(0,T;L^2(\Omega))$ and $u_0\in L^{2p}(\Omega)$, and moreover treating
a nonzero Dirichlet input without imposing any restriction on the space dimension
$d$. To the best of our knowledge, the smoothing effect and the analyticity exhibited in Propositions
\ref{p1} and \ref{p2} are the first results of this type stated in such a general context and, even for $F\equiv0$ and
$g\equiv0$, Theorem \ref{thm:main} is the first result of this type stated with an initial condition $u_0$ lying only in $L^{2p}(\Omega)$.
All existing results that we are aware of require at least that {$u_0\in H^s(\Omega)$
for some $s>\frac{d}{2}$} (see e.g.  \cite{HNWY,Yamamoto:2020order}).

The rest of the paper is organized as follows. In Section \ref{sec:prelim}, we present preliminary
results, i.e., analyticity and asymptotics of the solution $u$ to problem \eqref{eq1} for an admissible
tuple. The proof of Theorem \ref{thm:main} is given in Section \ref{sec:theorem}. Several numerical
tests are given in Section \ref{sec:numer} to illustrate the feasibility of unique order recovery. Throughout, the notation $C$ denotes a generic
positive constant independent of $t$ and it may change from line to line. Further, we often
write a bivariate function $f(x,t)$ as a vector valued function $f(t)$, by suppressing the dependence
on $x$. We denote by $L^2(\Omega;\rho \d x)$ the space of measurable functions $v$ satisfying
$\int_\Omega |v|^2\rho \d x<\infty$ endowed with the inner product $\langle u,v\rangle_{L^2(\Omega;\rho \d x)}
=\int_\Omega uv\rho \d x.$ Note that under condition \eqref{eqn:rho}, we have $L^2(\Omega;\rho \d x)
=L^2(\Omega)$ in the sense of set but equipped with different inner products and norms,
which are nonetheless equivalent to each other (under the given condition \eqref{eqn:rho} on $\rho$),
and thus we distinguish only the inner products but not the spaces.

\section{Preliminary properties}\label{sec:prelim}
In this section, we consider the direct problem \eqref{eq1} with an admissible tuple $(\alpha,\omega,a,q,
\rho ,u_0,F,g)$, and show the analyticity and asymptotic behavior of solutions $u$ of problem
\eqref{eq1} as $T\to+\infty$, using the standard separation of variables technique and Mittag-Leffler
functions as in \cite{SY}. These results will play a central role in the proof of Theorem
\ref{thm:main} in Section \ref{sec:theorem}.

\subsection{Mittag-Leffler function}
We shall use extensively the two-parameter Mittag-Leffler function $E_{\alpha,\beta}(z)$ defined by
\cite{KilbasSrivastavaTrujillo:2006,P,Jin:2021book}
\begin{equation*}
  E_{\alpha,\beta}(z) = \sum_{n=0}^\infty \frac{z^n}{\Gamma(n\alpha+\beta)},\quad z\in \mathbb{C}.
\end{equation*}
This function generalizes the exponential function $e^z$ in that $E_{1,1}(z)=e^z$, and it is an entire function of order $\frac1\alpha$ and type 1.
It has the following important asymptotic decay behavior in a sector of the complex plane $\mathbb{C}$
containing the negative real axis; (see \cite[pp. 34--35]{P} or \cite[Section 3.1]{Jin:2021book} for the proof).
\begin{lemma}\label{lem:mlf}
Let $\alpha\in(0,2)$, $\beta\in \mathbb{R}$, and $\mu\in (\frac\alpha2\pi,\min(\pi,\alpha\pi))$. Then
for any $ \mu\leq |\arg(z)|\leq \pi$ and $p\in \mathbb{N}$, there hold
\begin{align*}
 |E_{\alpha,\beta}(z)| &\leq c(1+|z|)^{-1} ,\\
 E_{\alpha,\beta}(z) & = -\sum_{k=1}^p\frac{z^{-k}}{\Gamma(\beta-k\alpha)} + \mathcal{O}(|z|^{-p-1}),\quad \mbox{as }|z|\to\infty.
\end{align*}
\end{lemma}
In Lemma \ref{lem:mlf} and below, since the set $\mathbb Z\setminus\mathbb N$
corresponds to the set of poles of the meromorphic extension to $\mathbb C$ of the
Gamma function $\Gamma(z)$, we use the convention $\frac{1}{\Gamma(m)}=0$, $m\in\mathbb Z\setminus\mathbb N.$

\subsection{Analyticity of solutions of problem \eqref{eq1}}
Consider the operator $A=\rho^{-1}\mathcal A$ acting on the space $L^2(\Omega;\rho \d x)$ with its domain
$H^2(\Omega)\cap H_0^1(\Omega)$. Then for any $s>0$, we may define the fractional power $A^s$
by spectral decomposition. Let $(\varphi_n)_{n\geq1}$ be an $L^2(\Omega;\rho\d x)$ orthonormal
basis of eigenfunctions of the operator $A$ associated with the non-decreasing sequence of eigenvalues
$(\lambda_n)_{n\geq1}$ (with multiplicity counted) of $A$. Then the operator $A^s$ is defined by
\begin{equation*}
  A^s v = \sum_{n=1}^\infty \lambda_n^s \langle v,\varphi_n\rangle_{L^2(\Omega;\rho \d x)}\varphi_n,
  \quad \mbox{with }  D(A^s) = \Big\{v\in L^2(\Omega):\sum_{n=1}^\infty \lambda_n^{2s}\langle v,\varphi_n\rangle_{L^2(\Omega;\rho \d x)}^2<\infty\Big\},
\end{equation*}
and the associated graph norm
$\|v\|_{D(A^s)} = (\sum_{n=1}^\infty \lambda_n^{2s}\langle v,\varphi_n\rangle_{L^2(\Omega;\rho \d x)}^2)^\frac12.$

Then we have the following result on the analytic extension of the solution $u$.
\begin{proposition}\label{p1}
Let $(\alpha,\omega,a,q,\rho,u_0,F,g)$ be an admissible tuple with $g\equiv0$.
Then the solution $u$ of problem \eqref{eq1} can be extended to a map $\tilde{u}\in L^1_{loc}(0,+\infty;L^2(\Omega))$ whose restriction to $\Omega\times (T_1,+\infty)$ is analytic with respect to $t\in(T_1,+\infty)$
as a function taking values in $C^1(\overline{\Omega})$.
\end{proposition}
\begin{proof}
First, for $t\in(0,+\infty)$, we define the maps $u_1$ and $u_2$ by
\begin{equation}\label{ppp1a}\begin{aligned}
 u_1(t)&=\sum_{n=1}^\infty \int_0^{\min(t,T)}(t-s)^{\alpha-1}E_{\alpha,\alpha}(-\lambda_n(t-s)^{\alpha})\langle F(s),\varphi_{n}\rangle_{L^2(\Omega;\rho \d x)}\d s\varphi_{n},\\
 u_2(t)&=\sum_{n=1}^\infty E_{\alpha,1}(-\lambda_nt^{\alpha})\langle u_0,\varphi_{n}\rangle_{L^2(\Omega;\rho \d x)}\varphi_n.
\end{aligned}\end{equation}
Let $\tilde{u}=u_1+u_2$. Then $\tilde{u}\in L^1_{loc}(0,+\infty;L^2(\Omega))$.
Moreover, according to \cite[Theorem 2.3]{KY3}, the map $\tilde{u}$ extends the solution $u$ to
problem \eqref{eq1}. It remains to show that the restriction of $u_j$, $j=1,2$, to $\Omega\times
(T_1,+\infty)$ is analytic with respect to $t\in(T_1,+\infty)$ as a function taking values in
$C^1(\overline{\Omega})$. We fix $\epsilon>0$ arbitrarily chosen and prove that $u_j$, $j=1,2$, is
analytic with respect to $t\in(T_1+\epsilon,+\infty)$ as a function taking values in $C^1
(\overline{\Omega})$. Since this result can be easily deduced for $\alpha=1$, we consider only
the case $\alpha\in(0,2)\setminus\{1\}$. Under condition \eqref{eqn:source}, we have
\begin{equation}\label{p1a}
u_1(t)=\sum_{n=1}^\infty \int_0^{T_1}(t-s)^{\alpha-1}E_{\alpha,\alpha}(-\lambda_n(t-s)^{\alpha})\left\langle F(s),\varphi_n\right\rangle_{L^2(\Omega;\rho \d x)}\d s\varphi_{n},\quad t\in(T_1,+\infty).
\end{equation}
Fix   $\ell_1=\left\lceil\frac{d}{4}\right\rceil$, $\theta\in(0,\min(\frac{(2-\alpha)\pi}{4\alpha},\frac{\pi}{4}))$, and $\mathcal D_\theta=\{T_1+\epsilon+re^{{\rm i}\beta}:\ \beta\in(-\theta,\theta),\ r>0\}.$
By Lemma \ref{lem:mlf}, with $\mu\in (\frac\alpha2\pi,\min(\pi,\alpha\pi))$, for all $z_1\in{\mathbb{S}_{\mu,\pi}:=\{z\in \mathbb C: \mu<|\textrm{arg}(z)|<\pi\}}$, we have
\begin{align*}
 E_{\alpha,\alpha}(z_1) &= -\sum_{k=1}^{\ell_1+1}\frac{z_1^{-k}}{\Gamma((1-k)\alpha)} + \mathcal{O}(|z_1|^{-\ell_1-2})\\
   &= -\sum_{k=1}^{\ell_1}\frac{z_1^{-1-k}}{\Gamma(-k\alpha)} + \mathcal{O}(|z_1|^{-\ell_1-2})\quad \mbox{as }|z_1|\to\infty,
\end{align*}
since $\frac{1}{\Gamma(0)}=0$.
Therefore, for fixed $\delta>0$, for all $z_1\in\mathbb S_{\mu,\pi}$ with $|z_1|\geq\delta$, we obtain
\begin{equation}\label{add}
\left|E_{\alpha,\alpha}(z_1)+\sum_{k=1}^{\ell_1}\frac{z_1^{-1-k}}{\Gamma(-k\alpha)}\right|\leq C|z_1|^{-\ell_1-2},
\end{equation}
with $C>0$ independent of $z_1$.
Note that one can find $\theta_0\in(0,\min(\frac{(2-\alpha)\pi}{4\alpha},\frac{\pi}{4}))$ such that for all
$s\in(0,T_1)$, $z\in \mathcal D_{\theta_0}$ and $n\in\mathbb N$, we have $-\lambda_n(z-s)^\alpha\in \mathbb{S}_{\mu,\pi}$ and $-\lambda_nz^\alpha\in \mathbb{S}_{\mu,\pi}$. In addition, for all $r>0$ and all $\beta\in(-\theta_0,\theta_0)$, we have
$|T_1+\epsilon+re^{{\rm i}\beta}|\geq |T_1+\epsilon+r\cos\beta|$
and since $0<\theta_0<\frac{\pi}{4}$, we deduce
$|T_1+\epsilon+re^{{\rm i}\beta}|\geq T_1+\epsilon+r\cos\theta_0\geq T_1+\epsilon.$
It follows that for all $z\in \mathcal D_{\theta_0}$, we have $|z|\geq T_1+\epsilon$, and thus,
$$|-\lambda_n(z-s)^\alpha|=\lambda_n|z-s|^\alpha\geq \lambda_n(|z|-s)^\alpha\geq\lambda_1\epsilon^\alpha>0,\quad z\in \mathcal D_{\theta_0},\ s\in(0,T_1),\ n\in\mathbb N.$$
Therefore, for all
$s\in(0,T_1)$, $z\in \mathcal D_{\theta_0}$ and $n\in\mathbb N$, we can apply \eqref{add} with $z_1=-\lambda_n(z-s)^\alpha$ and deduce
with $C>0$ independent of $s\in(0,T_1)$, $z\in \mathcal D_{\theta_0}$ and $n\in\mathbb N$,
$$\left|E_{\alpha,\alpha}(-\lambda_n(z-s)^\alpha)+\sum_{k=1}^{\ell_1}\frac{[-\lambda_n(z-s)^\alpha]^{-1-k}}{\Gamma(-k\alpha)}\right|\leq C|\lambda_n(z-s)^\alpha|^{-\ell_1-2}.$$
Multiplying both side of this inequality by $|z-s|^{\alpha-1}$, we obtain
\begin{align}\label{pppp6b}
&\left|(z-s)^{\alpha-1}E_{\alpha,\alpha}(-\lambda_n(z-s)^\alpha)+\sum_{k=1}^{\ell_1}\frac{(-1)^{k+1}(z-s)^{-k\alpha-1}}{\Gamma(-k\alpha)\lambda_n^{k+1}}\right|
\leq C\frac{|z-s|^{-(\ell_1+1)\alpha-1}}{\lambda_n^{\ell_1+2}},
\end{align}
 for all  $s\in(0,T_1)$, $z\in \mathcal D_{\theta_0}$ and $n\in\mathbb N$.
 Note that  for all $z\in \mathcal D_{\theta_0}$, we have $|z|\geq T_1+\epsilon$. Fixing $\gamma>0$, we deduce that, for all $z\in \mathcal D_{\theta_0}$ and $s\in(0,T_1)$,
$$|z-s|^{-\gamma}\leq (|z|-s)^{-\gamma}\leq (|z|-T_1)^{-\gamma}= |z|^{-\gamma}\left(1-\tfrac{T_1}{|z|}\right)^{-\gamma}\leq \left(1-\tfrac{T_1}{T_1+\epsilon}\right)^{-\gamma}|z|^{-\gamma}.$$
Combining this with the  estimate \eqref{pppp6b} gives
\begin{equation}\label{p6b}
\abs{(z-s)^{\alpha-1}E_{\alpha,\alpha}(-\lambda_n(z-s)^{\alpha})+\sum_{k=1}^{\ell_1}\frac{(-1)^{k+1}(z-s)^{-k\alpha-1}}{\Gamma(-k\alpha)\lambda_n^{k+1}}}\leq C\frac{|z|^{-1-(\ell_1+1)\alpha}}{\lambda_n^{\ell_1+2}}.
\end{equation}
Similarly, by Lemma \ref{lem:mlf}, for all $z_1\in\mathbb S_{\mu,\pi}$ with $|z_1|\geq\delta$, we obtain
$$\left|E_{\alpha,1}(z_1)+\sum_{k=1}^{\ell_1}\frac{z_1^{-k}}{\Gamma(1-k\alpha)}\right|\leq C|z_1|^{-\ell_1-1},$$
with $C>0$ independent of $z_1$. For all $z\in D_{\theta_0}$ and $n\in\mathbb N$, since $|z|\geq T_1+\epsilon$, we deduce
$$|-\lambda_nz^{\alpha}|= \lambda_n|z|^\alpha\geq \lambda_1(T_1+\epsilon)^\alpha>0.$$
Therefore, choosing $z_1=-\lambda_nz^{\alpha}\in \mathbb{S}_{\mu,\pi}$ and applying the above estimate lead to
\begin{equation}\label{p6bb}
\abs{E_{\alpha,1}(-\lambda_nz^{\alpha})+\sum_{k=1}^{\ell_1}\frac{(-1)^kz^{-k\alpha}}{\Gamma(1-k\alpha)\lambda_n^k}} \leq C\frac{|z|^{-(1+\ell_1)\alpha}}{\lambda_n^{\ell_1+1}}.
\end{equation}
For all $z\in \mathcal D_{\theta_0}$ and all $n\in\mathbb N$, let
\begin{align*}
X_n(z)&=\int_0^{T_1}(z-s)^{\alpha-1} E_{\alpha,\alpha}(-\lambda_n(z-s)^{\alpha})\langle F(s),\varphi_{n}\rangle_{L^2(\Omega;\rho \d x)}\d s\\
 &\qquad +\sum_{k=1}^{\ell_1}\int_0^{T_1}\frac{{(-1)^{k+1}}(z-s)^{-k\alpha-1}\left\langle v_k(s),\varphi_n\right\rangle_{L^2(\Omega;\rho \d x)}}{\Gamma(-k\alpha)}\d s,\\
Y_n(z)&=E_{\alpha,1}(-\lambda_nz^{\alpha})\langle u_0,\varphi_{n}\rangle_{L^2(\Omega;\rho \d x)}+\sum_{k=1}^{\ell_1}\frac{{(-1)^{k}}z^{-k\alpha}\left\langle w_k,\varphi_n\right\rangle_{L^2(\Omega;\rho \d x)}}{\Gamma(1-k\alpha)},
\end{align*}
with $v_k=A^{-k-1}F$ and $w_k=A^{-k}u_0$, for $k=1,\ldots,\ell_1$.
One can check that for all $n\in\mathbb N$, the maps $X_n$ and $Y_n$  are holomorphic
on $\mathcal D_{\theta_0}$. Moreover, for all $t>T_1+\epsilon$, we get
\begin{align*}
  u_1(t)+\sum_{k=1}^{\ell_1}\int_0^{T_1}\frac{{(-1)^{k+1}}(t-s)^{-k\alpha-1}v_k(s)}{\Gamma(-k\alpha)}\d s&=\sum_{n=1}^\infty X_n(t)\varphi_n,\\
  u_2(t)+\sum_{k=1}^{\ell_1}\frac{(-1)^{k}t^{-k\alpha}w_k}{\Gamma(1-k\alpha)}&=\sum_{n=1}^\infty Y_n(t)\varphi_n.
\end{align*}
Under the regularity assumptions on $\Omega$, $a$ and $q$ (from the admissible tuple), the space
$D(A^{\ell_1+1})$ continuously embeds into $H^{2\ell_1+2}(\Omega)$
\cite[Theorem 2.5.1.1]{Gr} and by Sobolev embedding theorem \cite{AdamsFournier:2003}, the space $D(A^{\ell_1+1})$
embeds continuously into $ C^1(\overline{\Omega})$. Therefore, applying \eqref{p6b}-\eqref{p6bb}, we deduce that,
for all $M,N\in \mathbb N$ and all $z\in \mathcal D_{\theta_0}$,
\begin{align}
   & \norm{\sum_{n=M}^N X_n(z)\varphi_n}_{ C^1(\overline{\Omega})}\leq C\norm{\sum_{n=M}^N X_n(z)\varphi_n}_{D(A^{\ell_1+1})}\nonumber\\
	\leq &C|z|^{-1-(\ell_1+1)\alpha}\norm{A^{-\ell_1-2} \left(\sum_{n=M}^N\left\langle F,\varphi_n\right\rangle_{L^2(\Omega;\rho \d x)}\varphi_n\right)}_{L^1(0,T;D(A^{\ell_1+1}))}\label{pp1b}\\
\leq &C|z|^{-1-(\ell_1+1)\alpha}\norm{\sum_{n=M}^N\left\langle F,\varphi_n\right\rangle_{L^2(\Omega;\rho \d x)}\varphi_n}_{ L^1(0,T;L^2(\Omega;\rho \d x))},\nonumber\\
   &\norm{\sum_{n=M}^N Y_n(z)\varphi_n}_{ C^1(\overline{\Omega})}\leq C\norm{\sum_{n=M}^N Y_n(z)\varphi_n}_{D(A^{\ell_1+1})}\nonumber\\
	\leq& C|z|^{-(1+\ell_1)\alpha}\norm{A^{-\ell_1-1}\left(\sum_{n=M}^N\left\langle u_0,\varphi_n\right\rangle_{L^2(\Omega;\rho \d x)}\varphi_n\right)}_{D(A^{\ell_1+1})}\label{pp1c}\\
\leq & C|z|^{-(1+\ell_1)\alpha}\norm{\sum_{n=M}^N\left\langle u_0,\varphi_n\right\rangle_{L^2(\Omega;\rho \d x)}\varphi_n}_{ L^2(\Omega;\rho \d x)},\nonumber
\end{align}
with $C>0$ being a constant independent of $M$, $N$ and $z$. The estimates \eqref{pp1b}-\eqref{pp1c} imply
that, for  any compact set $K\subset \mathcal D_{\theta_0}$, the sequences
$\sum_{n=1}^NX_n(z)\varphi_n$, $\sum_{n=1}^NY_n(z)\varphi_n$, for $N\in\mathbb N$,
converge uniformly with respect to $z\in K$ as  functions taking values in $C^1(\overline{\Omega})$. This proves that
the map $u_j^*$, $j=1,2$, given by
$u_1^*(t)=\sum_{n=1}^\infty X_n(t)\varphi_n$, $u_2^*(t)=\sum_{n=1}^\infty Y_n(t)\varphi_n$, for $t\in(T_1+\epsilon,+\infty)$,
are analytic as functions taking values in $C^1(\overline{\Omega})$. In addition,
since $u_0\in L^{2p}(\Omega)$ and $F\in L^1(0,T;L^p(\Omega))$, we deduce $v_k\in L^1(0,T;W^{4,p}(\Omega))$
and $w_k\in W^{2,2p}(\Omega)$,  $k=1,\ldots,\ell_1$ \cite[Theorem 2.5.1.1]{Gr}. This, the condition $p>\frac{d}{2}$ and
Sobolev embedding theorem give $v_k\in L^1(0,T;C^1(\overline{\Omega}))$ and $w_k\in C^1(\overline{\Omega})$,  $k=1,\ldots,\ell_1$.
Therefore, the maps
\begin{equation*}
 z\mapsto-\sum_{k=1}^{\ell_1}\int_0^{T_1}\frac{{(-1)^{k+1}}(z-s)^{-k\alpha-1}v_k(s)}{\Gamma(-k\alpha)}\d s,\quad z\mapsto-\sum_{k=1}^{\ell_1}\frac{{(-1)^{k}}z^{-k\alpha}w_k}{\Gamma(1-k\alpha)\lambda_n^k}
\end{equation*}
 are respectively  holomorphic extensions to $\mathcal D_{\theta_0}$ of the maps $u_1-u_1^*$ and $u_2-u_2^*$ restricted to $t\in(T_1+\epsilon,+\infty)$ as  functions taking values in $C^1(\overline{\Omega})$. Thus, both $u_1$ and $u_2$ are analytic with respect to  $t\in(T_1+\epsilon,+\infty)$ as  functions taking values in $C^1(\overline{\Omega})$. This proves that $\tilde{u}$ is analytic with respect to $t\in(T_1+\epsilon,+\infty)$ as a function taking values in $C^1(\overline{\Omega})$.
\end{proof}

We obtain a similar result for $F\equiv0$, $u_0\equiv0$ but $g\not\equiv0$.
\begin{proposition}\label{p2}
Let $(\alpha,\omega,a,q,\rho,u_0,F,g)$ be an admissible tuple, $u_0\equiv0$ and $F\equiv0$. Then the   solution
$u$ of problem \eqref{eq1} can be extended to a map $\tilde{u}\in L^1_{loc}(0,+\infty;L^2(\Omega))$ whose restriction to $\Omega\times (T_1,+\infty)$ is analytic with respect to $t\in(T_1,+\infty)$
as a function taking values in $C^1(\overline{\Omega})$.
\end{proposition}
\begin{proof}
Since the case for $\alpha=1$ is direct, we consider only the case
$\alpha\in(0,2)\setminus\{1\}$. We introduce the map for $t\in (0,+\infty)$,
\begin{equation}\label{uu}
\tilde{u}(t)=-\int_0^{{\min(t,T)}}(t-s)^{\alpha-1}\Big(\sum_{n=1}^\infty E_{\alpha,\alpha}(-\lambda_n(t-s)^{\alpha})\langle g(s),\partial_{\nu_a}\varphi_{n}\rangle_{L^2(\partial\widetilde{\Omega})}\varphi_{n}\Big)\d s,
\end{equation}
where $\langle\cdot,\cdot\rangle_{L^2(\partial\widetilde\Omega)}$ is the standard
$L^2(\partial\widetilde\Omega)$ inner product. Under the given condition on $g$, one can readily check that $\tilde{u}\in L^1_{loc}
(0,+\infty;L^2(\Omega))$ and, in view of \cite[Proposition 3.1]{KLLY}, there holds
$\tilde{u}=u$ on $\Omega\times(0,T)$. Thus, fixing $\epsilon>0$ arbitrarily chosen,
the proposition is proven if we show that $\tilde{u}$ is analytic with respect to $t\in(T_1+\epsilon,+\infty)$ as a function taking
values in $C^1(\overline{\Omega})$. Applying \eqref{eqn:diri-bc}, we find for $t\in (T_1,+\infty),$
\begin{equation}\label{p6a}
\tilde{u}(t)=-\int_0^{T_1}(t-s)^{\alpha-1}\Big(\sum_{n=1}^\infty E_{\alpha,\alpha}(-\lambda_n(t-s)^{\alpha})\langle g(s),\partial_{\nu_a}\varphi_{n}\rangle_{L^2(\partial\widetilde{\Omega})}\varphi_{n}\Big)\d s,
\end{equation}
where for $x\in \partial\widetilde \Omega$, the notation
$$\partial_{\nu_a }h(x) := \sum_{i,j=1}^d a_{i,j}(x) \partial_{x_j} h(x) \nu_i(x)$$
denotes the conormal derivative.
For each $t\in(0,T)$, let $G(\cdot,t)$ be the solution of
\begin{equation}\label{eq2}
\left\{ \begin{aligned}
\mathcal{A} G(\cdot,t) & =  0, && \mbox{in }\Omega,\\
 G(\cdot,t) & =   g(\cdot,t), && \mbox{on } \partial \widetilde{\Omega},\\
 G(\cdot,t) & =   0, && \mbox{on }  \partial \omega.
\end{aligned}
\right.
\end{equation}
Since $g\in  L^1(0,T; W^{2-\frac{1}{p},p}(\partial\widetilde{\Omega}))\cap_{k=0}^{\lceil\alpha\rceil} W^{\lceil\alpha\rceil-k,1}(0,T; H^{\frac{1}{2}+\frac{k}{\lceil\alpha\rceil}}(\partial\widetilde{\Omega}))$, by the standard
elliptic regularity theory, we have $G\in L^1(0,T;H^2(\Omega))\cap L^1(0,T;W^{2,p}(\Omega))$.
We fix also   $y_k(\cdot,t)=A^{-k}G(\cdot,t)$, $k=1,\ldots,\ell_1:=\left\lceil\frac{d}{4}\right\rceil$, and, by \cite[Lemma 2.1]{KY3}, we deduce for $k=1,\ldots,\ell_1$,
\begin{equation*}
   \left\langle y_k(t),\varphi_{n}\right\rangle_{L^2(\Omega;\rho \d x)}=\frac{\left\langle G(t),\varphi_{n}\right\rangle_{L^2(\Omega;\rho \d x)}}{\lambda_n^k}=-\frac{\left\langle g(t),\partial_{\nu_a}\varphi_{n}\right\rangle_{L^2(\partial\widetilde\Omega)}}{\lambda_n^{k+1}},\quad t\in(0,T),\ n\in\mathbb N.
\end{equation*}
The condition $G\in L^1(0,T;H^2(\Omega))$ implies that
the sequence
$$\sum_{n=1}^N\frac{\langle g(t),\partial_{\nu_a}\varphi_{n}\rangle_{L^2(\partial\widetilde\Omega)}}{\lambda_n}\varphi_n, \quad N\in\mathbb{N}, \,t\in(0,T),$$
converges in the sense of $L^1(0,T;L^2(\Omega))$.
For all $z\in \mathcal D_{\theta_0}$ and all $n\in\mathbb N$, let
\begin{align*}
H_n(z) &=-\int_0^{T_1}(z-s)^{\alpha-1} E_{\alpha,\alpha}(-\lambda_n(z-s)^{\alpha})\langle g(s),\partial_{\nu_a}\varphi_{n}\rangle_{L^2(\partial\widetilde{\Omega})}\d s\\
 &\quad +\sum_{k=1}^{\ell_1}\int_0^{T_1}\frac{{(-1)^{k+1}}(z-s)^{-k\alpha-1}\left\langle y_k(s),\varphi_n\right\rangle_{L^2(\Omega;\rho \d x)}}{\Gamma(-k\alpha)}\d s.
\end{align*}
Repeating the argument for Proposition \ref{p1}, for all $M,N\in \mathbb N$ and all $z\in \mathcal D_{\theta_0}$, we obtain
\begin{equation}\label{pr1c}
   \norm{\sum_{n=M}^N H_n(z)\varphi_n}_{ C^1(\overline{\Omega})}\leq C|z|^{-1-(\ell_1+1)\alpha}\norm{\sum_{n=M}^N\left\langle G,\varphi_n\right\rangle_{L^2(\Omega;\rho \d x)}\varphi_n}_{ L^1(0,T;L^2(\Omega;\rho \d x))},
\end{equation}
with $C>0$ being a constant independent of $M$, $N$ and $z$. Then, we deduce that the map
\begin{equation*}
u^*(t):=\sum_{n=1}^\infty H_n(t)\varphi_n,\quad t\in(T_1+\epsilon,+\infty)
\end{equation*}
is analytic as a function taking values in $C^1(\overline{\Omega})$. Similarly, since
$g\in L^1(0,T; W^{2-\frac{1}{p},p}(\partial\widetilde{\Omega}))\cap_{k=0}^{\lceil\alpha\rceil} W^{\lceil\alpha\rceil-k,1}(0,T; H^{\frac{1}{2}+\frac{k}{\lceil\alpha\rceil}}(\partial\widetilde{\Omega}))$,
we deduce $G\in L^1(0,T;W^{2,p}(\Omega))$ and, furthermore, applying \cite[Theorem 2.5.1.1]{Gr}, the condition $p>\frac{d}{2}$ and Sobolev embedding theorem \cite{AdamsFournier:2003}, we obtain $y_k\in L^1(0,T;W^{4,p}(\Omega))\hookrightarrow L^1(0,T;C^1(\overline{\Omega}))$, $k=1,\ldots,\ell_1$. Hence, repeating the argument for Proposition \ref{p2}, we deduce that $\tilde{u}$ is analytic with respect to $t\in(T_1+\epsilon,+\infty)$ as a function taking values in $C^1(\overline{\Omega})$.
\end{proof}

\subsection{Asymptotic properties of the analytic extension of solutions of problem \eqref{eq1}}
Now we consider the analytic extension $\tilde{u}\in C(T_1,+\infty; C^1(\overline{\Omega}))$ of the solution $u$ of \eqref{eq1} in Propositions \ref{p1} and \ref{p2}. Thus, for any  $x_0\in\partial\widetilde\Omega$, the map $(T_1,+\infty)\ni
t\mapsto \partial_\nu \tilde{u}(x_0,t)$ belongs to $ C(T_1,+\infty)$. Below we study the asymptotic behavior of $\partial_\nu \tilde{u}(x_0,t)$
as $t\to+\infty$, and analyze separately the three cases, i.e., $F\equiv0$ and $g\equiv0$, $u_0\equiv0$ and $g\equiv0$, and
$u_0\equiv0$ and $F\equiv0$. The next result gives the asymptotic, as $t\to+\infty$, for $t\mapsto\partial_\nu \tilde{u}(x_0,t)$
when $F\equiv 0$ and $g\equiv0$.

\begin{proposition}\label{prop:asympt-u0}
Let $(\alpha,\omega,a,q,\rho,u_0,F,g)$ be admissible tuple,
$F\equiv0$, $g\equiv0$ and $x_0\in\partial\widetilde\Omega$.
If $\partial_\nu A^{-1}u_0(x_0)\neq0$, then the extension $\tilde{u}$ of the solution $u$ of
problem \eqref{eq1} defined in  Proposition \ref{p1} satisfies
\begin{equation}\label{p2a}
\partial_\nu \tilde{u}(x_0,t)=-\frac{\partial_\nu A^{-1}u_0(x_0)}{\Gamma(1-\alpha)}t^{-\alpha}+{\mathcal O}(t^{-2\alpha}),\quad \mbox{as } t\to+\infty.
\end{equation}
\end{proposition}
\begin{proof}
Applying \eqref{pp1c} with $M=1$ and $N=\infty$, we deduce that, for all $t>T_1+1$,
$$\norm{\tilde{u}(t)+\sum_{k=1}^{\ell_1}\frac{{(-1)^{k}}t^{-k\alpha}w_k}{\Gamma(1-k\alpha)}}_{C^1(\overline{\Omega})}\leq Ct^{-\alpha(\ell_1+1)},$$
with $C>0$ independent of $t$ and $w_k=A^{-k}u_0$, $k=1,\ldots,\ell_1$. It follows that
$$\partial_\nu \tilde{u}(x_0,t)=-\sum_{k=1}^{\ell_1}\frac{{(-1)^{k}}t^{-k\alpha}\partial_\nu w_k(x_0)}{\Gamma(1-k\alpha)}+{\mathcal O}(t^{-\alpha(\ell_1+1)}),\quad \mbox{as } t\to+\infty$$
which clearly implies \eqref{p2a}.
\end{proof}

The next result gives the asymptotics of the map $t\mapsto\partial_\nu \tilde{u}(x_0,t)$ when $u_0\equiv0$ and $g\equiv0$.
\begin{proposition}\label{prop:asympt-f}
Let $(\alpha,\omega,a,q,\rho,u_0,F,g)$ be an admissible tuple, $u_0\equiv0$, $g\equiv0$
and $x_0\in\partial\widetilde\Omega$, and let $$F^*=\int_0^{T}F(t)\d t\in  L^p(\Omega)\cap L^2(\Omega).$$
If $\partial_\nu A^{-2}F^*(x_0)\ne0$, then the extension $\tilde{u}$ of the solution $u$ of
problem \eqref{eq1} defined in  Proposition \ref{p1} satisfies
\begin{equation}\label{p3a}
\partial_\nu \tilde{u}(x_0,t)=\frac{\partial_\nu A^{-2}F^*(x_0)}{\Gamma(-\alpha)}t^{-1-\alpha}+{\mathcal O}(t^{-1-2\alpha}),\quad\mbox{as }t\to+\infty.
\end{equation}
\end{proposition}
\begin{proof}
Applying \eqref{pp1b} with $M=1$ and $N=\infty$, we deduce that, for all $t>T_1+1$, we have
$$\norm{\tilde{u}(t)+\sum_{k=1}^{\ell_1}\int_0^{T_1}\frac{(-1)^{k+1}(t-s)^{-k\alpha-1}v_k(s)}{\Gamma(-k\alpha)}\d s}_{C^1(\overline{\Omega})}\leq Ct^{-1-(\ell_1+1)\alpha},$$
with $C>0$ independent of $t$ and $v_k=A^{-k-1}F$, $k=1,\ldots,\ell_1$. Combining this with the fact that, for all $k=1,\ldots,\ell_1$, $v_k\in L^1(0,T;\mathcal C^{1}(\overline{\Omega}))$, we obtain
\begin{equation*}
   \partial_\nu \tilde{u}(x_0,t)=\partial_\nu \left(\int_0^{T_1}\frac{(t-s)^{-\alpha-1}v_1(s)}{\Gamma(-\alpha)}\d s\right)+{\mathcal O}(t^{-1-2\alpha}),\quad \mbox{as } t\to+\infty.
\end{equation*}
Further, we have
$(t-s)^{-1-\alpha}=t^{-1-\alpha}+{\mathcal O}(t^{-1-2\alpha})$, for $s\in(0,T_1)$, as $t\to+\infty$,
and hence,
\begin{equation*}
  \partial_\nu \tilde{u}(x_0,t)=t^{-1-\alpha}\partial_\nu\left(\int_0^{T_1}\frac{ v_1(s)}{\Gamma(-\alpha)}\d s\right)+{\mathcal O}(t^{-1-2\alpha}),\quad\mbox{as }t\to+\infty.
\end{equation*}
Finally, applying \eqref{eqn:source} and noting
$$\int_0^{T_1} v_1(s)\d s=\int_0^{T_1}A^{-2}F(s)\d s=A^{-2}F_*,$$
we obtain \eqref{p3a}.
\end{proof}

Let $G^*=\int_0^TG(s)\d s$, with $G$ solving \eqref{eq2}. Then $G^*\in H^2(\Omega)\cap W^{2,p}(\Omega)$ is the unique solution of the boundary value problem
\begin{equation}\label{eqn:Gstar}
\left\{ \begin{aligned}
\mathcal{A}G^* & =  0,  &&\mbox{in } \Omega,\\
 G^* & =  g^*, && \mbox{on }  \partial \widetilde\Omega,\\
 G^* & = 0, && \mbox{on }\partial\omega.
\end{aligned}\right.
\end{equation}
Combining this with the arguments in Proposition \ref{prop:asympt-f}  and applying estimate \eqref{pr1c} give
the asymptotics, as $t\to+\infty$, of the map $t\mapsto\partial_\nu u(x_0,t)$,
when  $u_0\equiv0$ and $F\equiv0$.
\begin{proposition}\label{prop:asympt-g}
Let $(\alpha,\omega,a,q,\rho,u_0,F,g)$ be an admissible tuple, $u_0\equiv0$, $F\equiv0$
and $x_0\in\partial\tilde{\Omega}$, and $G^* \in H^2(\Omega)\cap W^{2,p}(\Omega)$ solve \eqref{eqn:Gstar}. If $\partial_\nu A^{-1}G^*(x_0)\neq0$, then $\tilde{u}$ satisfies
\begin{equation}\label{p5a}
\partial_\nu \tilde{u}(x_0,t)=\frac{\partial_\nu A^{-1}G^*(x_0)}{\Gamma(-\alpha)}t^{-1-\alpha}+{\mathcal O}(t^{-1-2\alpha}),\quad\mbox{as }t\to+\infty.
\end{equation}
\end{proposition}

\begin{remark}
The proofs of Propositions \ref{prop:asympt-u0}--\ref{prop:asympt-g} indicate that one can actually obtain more precise asymptotic expansions
including high-order terms. For example, for $u_0\in L^{2p}(\Omega)$, $F\equiv0$ and $g\equiv0$, there holds
\begin{equation*}
\partial_\nu \tilde{u}(x_0,t)=\sum_{k=1}^{\ell_1}\frac{(-1)^k\partial_\nu A^{-k}u_0(x_0)}{\Gamma(1-k\alpha)}t^{-k\alpha}+{\mathcal O}(t^{-(\ell_1+1)\alpha}),\quad \mbox{as } t\to+\infty.
\end{equation*}
Nonetheless, under the conditions of Theorem \ref{thm:main}, the leading term in the expansion does not
vanish, cf. Lemma \ref{lem:nonzero}, and suffices the proof of Theorem \ref{thm:main}.
\end{remark}

\section{Proof of Theorem \ref{thm:main}}\label{sec:theorem}
In this section, we give the proof of Theorem \ref{thm:main}. To this end, for $k=1,2$, we define the operators
corresponding to $A_k=\rho_k^{-1}\mathcal A_k$ acting on
$L^2(\Omega_k;\rho_k\d x)$ with  their domain $D(A_k)=H^1_0(\Omega_k)\cap
H^2(\Omega_k)$. Further, for $k=1,2$, let
\begin{equation*}
   v^k=A_k^{-1}u_0^k\quad \mbox{and}\quad w^k = A_k^{-2}F_k^*+A_k^{-1}G_k^*,
\end{equation*}
with $F_k^*=\int_0^TF_k(t)\d t$ (cf. Proposition \ref{prop:asympt-f}) and
$G_k^* $ is defined in \eqref{eqn:Gstar} with $g^*=g_k^*$ on the domain
$\Omega_k$ (cf. Propositions \ref{p2} and \ref{prop:asympt-g}).
First we give an auxiliary result on $v^k$ and $w^k$.
\begin{lemma}\label{lem:nonzero}
The following statements hold.
\begin{itemize}
  \item[{\rm(i)}] If condition {\rm(}i{\rm)} of Theorem \ref{thm:main} holds, then  $v^k\in C^1(\overline{\Omega_k})$,
and for any $x_0\in\partial\widetilde{\Omega}$, $\partial_\nu v^k(x_0)\neq0$.
  \item[{\rm(ii)}] If condition {\rm(}ii{\rm)} of Theorem \ref{thm:main} holds, then  $w^k\in C^1(\overline{\Omega_k})$,
and for any $x_0\in\partial\widetilde{\Omega}$, $\partial_\nu w^k(x_0)\neq0$.
\end{itemize}
\end{lemma}
\begin{proof}
We suppress the subscript $k$ in the proof. The regularity $v\in C^1(\overline{\Omega})$ and $w\in C^1(\overline{\Omega})$
follows directly from Sobolev embedding theorem \cite{AdamsFournier:2003} and the elliptic regularity property (see e.g.
\cite[Theorem 2.5.1.1]{Gr}). Under condition (i), $u_0$ is of constant sign, and we may assume that $u_0\leq0$.
Note that the function $v$ solves
$$\left\{ \begin{aligned}
\mathcal A v & =  \rho u_0, && \mbox{in }\Omega,\\
v(x) & =   0, &&\mbox{on }\partial \Omega.
\end{aligned}
\right.$$
Since  $\mathcal A v=\rho u_0\leq 0$ in $\Omega$, $u_0\not\equiv0$ and $v|_{\partial\Omega}=0$, the
strong maximum principle \cite[Theorem 3.5]{GT} implies
$v(x)<0=v(x_0)$, for $x\in\Omega,\ x_0\in\partial\widetilde{\Omega}\subset\partial \Omega$.
Thus, the Hopf lemma  \cite[Lemma 3.4]{GT} implies $\partial_\nu v(x_0)>0$ for $x_0\in
\partial\widetilde{\Omega}.$ This shows the assertion in (i).

\noindent Now we turn to condition (ii).  Since $F^*$ and $g^*$ have the same constant sign, we may assume that $F^*\leq0$ and $g^*\leq0$. Let $y=A^{-1}F^*$. Then the function
$w$ solves
$$\left\{ \begin{aligned}
\mathcal A w & =  \rho y+\rho G^*, && \mbox{in } \Omega,\\
w & =   0, &&\mbox{on }\partial \Omega.
\end{aligned}
\right.$$
Since $F^*\in L^p(\Omega)$, by \cite[Theorem 2.4.2.5]{Gr}, there holds $y\in W^{2,p}(\Omega)$ and, since $p>\frac{d}{2}$
by assumption, the Sobolev embedding theorem implies that $y\in C(\overline{\Omega})$.  This, the fact $F^*\leq0$, $F^*\not\equiv0$ and the maximum principle
\cite[Corollary 3.2]{GT} imply $y\leq0$. Similarly, we can prove $G^*\leq0$ and it follows $\max(\rho y,\rho G^*)\leq0$. Moreover, the fact that $F^*\not\equiv0$ or $g^*\not\equiv0$ implies that $\rho y\not\equiv0$ or $\rho G^*\not\equiv0$. Thus $\rho y+\rho G^*\leq0$ and $\rho y+\rho G^*\not\equiv0$. Consequently, by repeating the above application of the strong maximum principle and the Hopf lemma, we deduce that, for all $x_0\in\partial\tilde{\Omega}$, we  have $\partial_\nu w(x_0)>0$.
\end{proof}

Now we can give the proof of Theorem \ref{thm:main}.
\begin{proof}
Let $u^k$, $k=1,2$, be the extension, introduced in Proposition \ref{p1} and \ref{p2}, of the solution of problem \eqref{eq1}
corresponding to the admissible tuple $(\alpha_k,\omega_k,a^k,\rho_k,q_k,u_0^k,F_k,g_k)$.
For all $\delta\in(0,T-T_1)$, the regularity $u^k\in C([T-\delta,T];C^1(\overline{\Omega_k}))$ is direct from
Propositions \ref{p1} and \ref{p2}. Thus it suffices to show the uniqueness. Fix $\delta\in(0,T-T_1)$, $x_0\in\partial\widetilde{\Omega}$
and let condition \eqref{t1a} be fulfilled. From Propositions \ref{p1} and \ref{p2}, we
deduce that $(T_1,+\infty)\ni t\mapsto \partial_\nu u^k(x_0,t)$, $k=1,2$, is an analytic function. {Moreover, following the discussions at the beginning of Proposition \ref{p1}, one can check that the restriction of $u^k$, $k=1,2$, to $\Omega\times (0,T)$ coincides with the solution of \eqref{eq1} corresponding to the admissible tuple $(\alpha_k,\omega_k,a^k,\rho_k,q_k,u_0^k,F_k,g_k)$.}
Therefore, condition \eqref{t1a} and unique continuation of analytic functions imply
\begin{equation}\label{t1b}
   \partial_\nu u^1(x_0,t)=\partial_\nu u^2(x_0,t),\quad t\in(T-\delta,+\infty).
\end{equation}
It remains to show that the identity \eqref{t1b} and one of the conditions (i) and (ii)
imply $\alpha_1=\alpha_2$.
First, we prove Theorem \ref{thm:main} under condition (i). For $k=1,2$, let $(\varphi_n^k)_{n\geq1}$ be an $L^2(\Omega;\rho_k\d x)$ orthonormal
basis of eigenfunctions of the operator $A_k$ associated with the non-decreasing sequence of eigenvalues
$(\lambda_n^k)_{n\geq1}$. We recall  that
$u^k=y_1^k+y_2^k+y_3^k$, with
\begin{align*}
y_1^k(t)&=\sum_{n=1}^\infty E_{\alpha,1}(-\lambda_n^kt^{\alpha})\langle u_0^k,\varphi_{n}\rangle_{L^2(\Omega_k;\rho_k \d x)}\varphi_n^k,\\
y_2^k(t)&=\sum_{n=1}^\infty \int_0^{\min(t,T)}(t-s)^{\alpha-1}E_{\alpha,\alpha}(-\lambda_n^k(t-s)^{\alpha})\langle F_k(s),\varphi_{n}^k\rangle_{L^2(\Omega_k;\rho_k \d x)}\d s\varphi_{n}^k,\\
y_3^k(t)&=-\int_0^{{\min(t,T)}}(t-s)^{\alpha-1}\left(\sum_{n=1}^\infty E_{\alpha,\alpha}(-\lambda_n^k(t-s)^{\alpha})\langle g_k(s),\partial_{\nu_{a_k}}\varphi_{n}^k\rangle_{L^2(\partial\widetilde{\Omega})}\varphi_{n}^k\right)\d s,
\end{align*}
for $t\in (0,+\infty)$. Propositions \ref{prop:asympt-u0}, \ref{prop:asympt-f} and \ref{prop:asympt-g} yield that for $t\to{+}\infty$
\begin{align*}
   \partial_\nu y_1^k(x_0,t)&=-\frac{\partial_\nu v^k(x_0)}{\Gamma(1-\alpha_k)}t^{-\alpha_k}+{\mathcal O}(t^{-2\alpha_k}),\\
   \partial_\nu y_2^k(x_0,t)+\partial_\nu y_3^k(x_0,t)&=\frac{\partial_\nu w^k(x_0)}{\Gamma(-\alpha_k)}t^{-1-\alpha_k}+{\mathcal O}(t^{-1-2\alpha_k}).
\end{align*}
Therefore, we find for $t\to{+}\infty$,
\begin{equation*}
   \partial_\nu u^k(x_0,t)=-\frac{\partial_\nu v^k(x_0)}{\Gamma(1-\alpha_k)}t^{-\alpha_k}+{\mathcal O}(t^{-\min(2\alpha_k,1+\alpha_k)}),
\end{equation*}
and the identity \eqref{t1b} implies
\begin{equation}
\label{t1d}
-\frac{\partial_\nu v^1(x_0)}{\Gamma(1-\alpha_1)}t^{-\alpha_1}+{\mathcal O}(t^{-\min(2\alpha_1,1+\alpha_1)})=-\frac{\partial_\nu v^2(x_0)}{\Gamma(1-\alpha_2)}t^{-\alpha_2}+{\mathcal O}(t^{-\min(2\alpha_2,1+\alpha_2)}).
\end{equation}
Combining this with the fact that $\partial_\nu v^k(x_0)\neq0$, $k=1,2$, cf. Lemma \ref{lem:nonzero}(i),
one can easily prove by contradiction that \eqref{t1d} implies $\alpha_1=\alpha_2$. Note that here for
$\alpha_1=1$ or $\alpha_2=1$ one can deduce $\alpha_1=1=\alpha_2$, since
 $\mathbb{Z}\setminus\mathbb{N}$ is the set of poles of $\Gamma(z)$ in \eqref{t1d}.

\noindent Next, assume that condition (ii) holds. The preceding argumentation gives that for $t\to+\infty$
\begin{equation*}
   \partial_\nu u^k(x_0,t)=\frac{\partial_\nu w^k(x_0)}{\Gamma(-\alpha_k)}t^{-1-\alpha_k}+{\mathcal O}(t^{-1-2\alpha_k}).
\end{equation*}
Then, the condition \eqref{t1b} implies
\begin{equation}\label{t1e}
\frac{\partial_\nu w^1(x_0)}{\Gamma(-\alpha_1)}t^{-1-\alpha_1}+{\mathcal O}(t^{-1-2\alpha_1})=\frac{\partial_\nu w^2(x_0)}{\Gamma(-\alpha_2)}t^{-1-\alpha_2}+{\mathcal O}(t^{-1-2\alpha_2}),\quad\mbox{as }t\to+\infty.
\end{equation}
This, the fact that $\partial_\nu w^k(x_0)\neq0$, $k=1,2$, cf. Lemma \ref{lem:nonzero}(ii), and \eqref{t1e}
imply that $\alpha_1=\alpha_2$.
\end{proof}

\begin{remark}
If the inclusion $\omega=\emptyset$, the results in Theorem \ref{thm:main} hold also for $d=1$.
\end{remark}

\section{Numerical experiments and discussions}\label{sec:numer}
Now we discuss the numerical recovery of the fractional order $\alpha$ from the flux data $\partial_\nu u(x_0,t)$
over the observation window $[T_1,T_2]$,
which has not been extensively studied in the literature so far. Hatano et al \cite{HatanoHatano:1998}
employed the asymptotic formula and numerical differentiation to recover the order $\alpha$. We describe a numerical
procedure motivated by the analysis in Section \ref{sec:theorem}. The analysis in Section \ref{sec:theorem}
proceeds in two steps: (i) analytic continuation and (ii) asymptotic matching. The first step can be numerically
ill-conditioned, especially when the measurement time horizon $[T_1,T_2]$ is small or $T_1$ is very small.
Nonetheless, when the observation time $T_1$ is sufficiently large, there is a simple recipe to recover the
fractional order $\alpha$. Specifically, for large $t$, when $u_0\not\equiv 0$, the normal derivative
$\partial_\nu u(x_0,t)$ behaves like
\begin{equation}\label{eqn:exp}
 {h(t)}\equiv \partial_\nu u(x_0,t) = c_1t^{-\alpha} + c_2t^{-2\alpha} + c_3 t^{-1-\alpha} + c_4t^{-1-2\alpha} + h.o.t.
\end{equation}
Thus, one may recover the order $\alpha$ by fitting to a mixture of powers $\{t^{-k\alpha},
t^{-1-k\alpha}\}_{k=1}^\infty$. This can be done with the following nonlinear least-squares problem
\begin{equation}\label{eqn:lsq}
  (\alpha^*,\mathbf{c}^*) = \arg\min_{\alpha\in[0,2],\mathbb{c}}\sum_{i=1}^N\Big({h}(t_i)-\sum_{k=1}^Kc_kt^{-\alpha_k}\Big)^2,
\end{equation}
with $\alpha_k=-k\alpha$ or $\alpha_k=1-k\alpha$, depending on the \textit{a priori} knowledge on $u_0$ (i.e.,
condition (i) or (ii) / (iii) in Theorem \ref{thm:main}), $\mathbf{c}=(c_1,\ldots,c_K)^T \in \mathbb{R}^K$,
and $\{t_i\}_{i=1}^N$ are the sampling points at which discrete observations are available. The formulation
\eqref{eqn:lsq} is very flexible, and capable of handling sparsely / irregularly sparse data points. Note that we do not
include a penalty term (e.g., $\ell^2$ or $\ell^1$) in the formulation \eqref{eqn:lsq}, since we generally take only a few terms in the
expansion \eqref{eqn:exp}, which has a built-in regularizing effect. {In addition, the optimal strength
of the penalty should depend on the noise magnitude, which differs dramatically for different cases (e.g.,
the presence of a nonzero initial condition). Hence the use of a penalty requires much tuning in the current
context, and we do not pursue the penalized approach in this work.} The optimization problem in \eqref{eqn:lsq}
can be readily solved by any stand-alone optimizer, e.g., limited-memory BFGS.
Note that the exponent $\alpha$ can be warm started by estimating with one single term (for which the problem
can be solved explicitly with log transformation), which can often deliver reasonable estimates. Numerically,
we observe that the procedure is fairly robust.

Below we present several numerical tests to show the feasibility of the approach. In all the experiments
below, the density $\rho$ is fixed at $\rho\equiv1$. The exact flux data ${h}(t)=\partial_\nu u(x_0,t)$ is generated by solving the direct
problem \eqref{eq1} over a large time interval $[0,T]$ with $T\gg 1$, which is fixed at $T=100$ below,
discretized with a time step size $\tau = 1\times 10^{-4}$. The simulation of the direct problem requires extra care in the practical
implementation in order to be numerically efficient, since the straightforward implementation of time stepping schemes incur
huge time and storage issues. We employ the well known sum of exponentials approximation of the singular kernel to speed up the computation;
see the appendix for details. The order $\alpha$ is recovered using eleven discrete observations
that are equally spaced within the window $[T_1,T_2]$. The noisy data ${h}^\delta $ is generated by
adding componentwise noise to the exact data ${h}(t)$ by ${h}^\delta(t_i)={h}(t_i) (1+\epsilon\xi(t_i))$, where
$\epsilon$ denotes the relative noise level and the noise $\xi(t_i)$ follows the standard Gaussian distribution.
Since the subdiffusion and diffusion wave cases exhibit different behavior, we discuss the corresponding
numerical results separately.

\subsection{Numerical results for subdiffusion}

First we present one-dimensional examples, one with nonzero initial condition, and the
other two with zero initial condition. The notation $\chi_S$ denotes the characteristic function
of a set $S$.
\begin{example}\label{exam:1dsub}
The domain $\Omega$ is taken to be the unit interval $[0,1]$, and $\omega=\emptyset$. The
observation point $x_0$ is the left end point $x_0=0$.
\begin{itemize}
  \item[{\rm(i)}] $(a,q,u_0,F,g)=(1+x^2,1,x^2(1-x),e^{x(1-x)}x(1-x)t\chi_{[0,0.1]}(t),0)$.
  \item[{\rm(ii)}] $(a,q,u_0,F,g)=(1,1+\sin(x),0,e^{x^2}\sin(\pi x)\chi_{[0,0.1]}(t),0)$.
  \item[{\rm(iii)}] $(a,q,u_0,F,g)=(1+\sin(\pi x),\cos(\pi x),0,0,e^{t}\chi_{[0,0.1]}(t))$, with
  the Dirichlet input $g$ specified on the left end point $x=0$.
\end{itemize}
\end{example}

\begin{figure}[hbt!]
\centering
\setlength{\tabcolsep}{0pt}
\begin{tabular}{ccc}
  \includegraphics[width=0.33\textwidth]{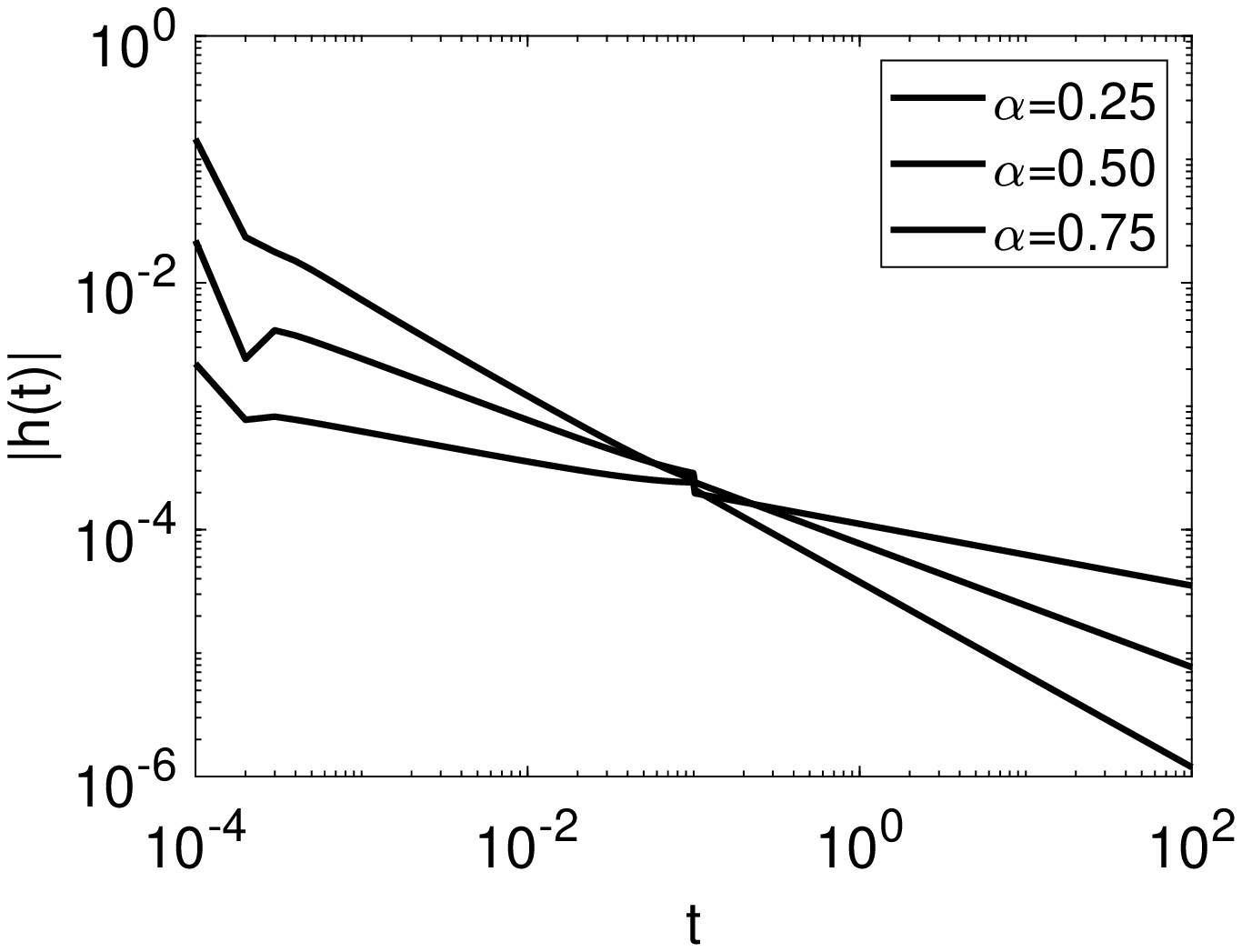} &
  \includegraphics[width=0.33\textwidth]{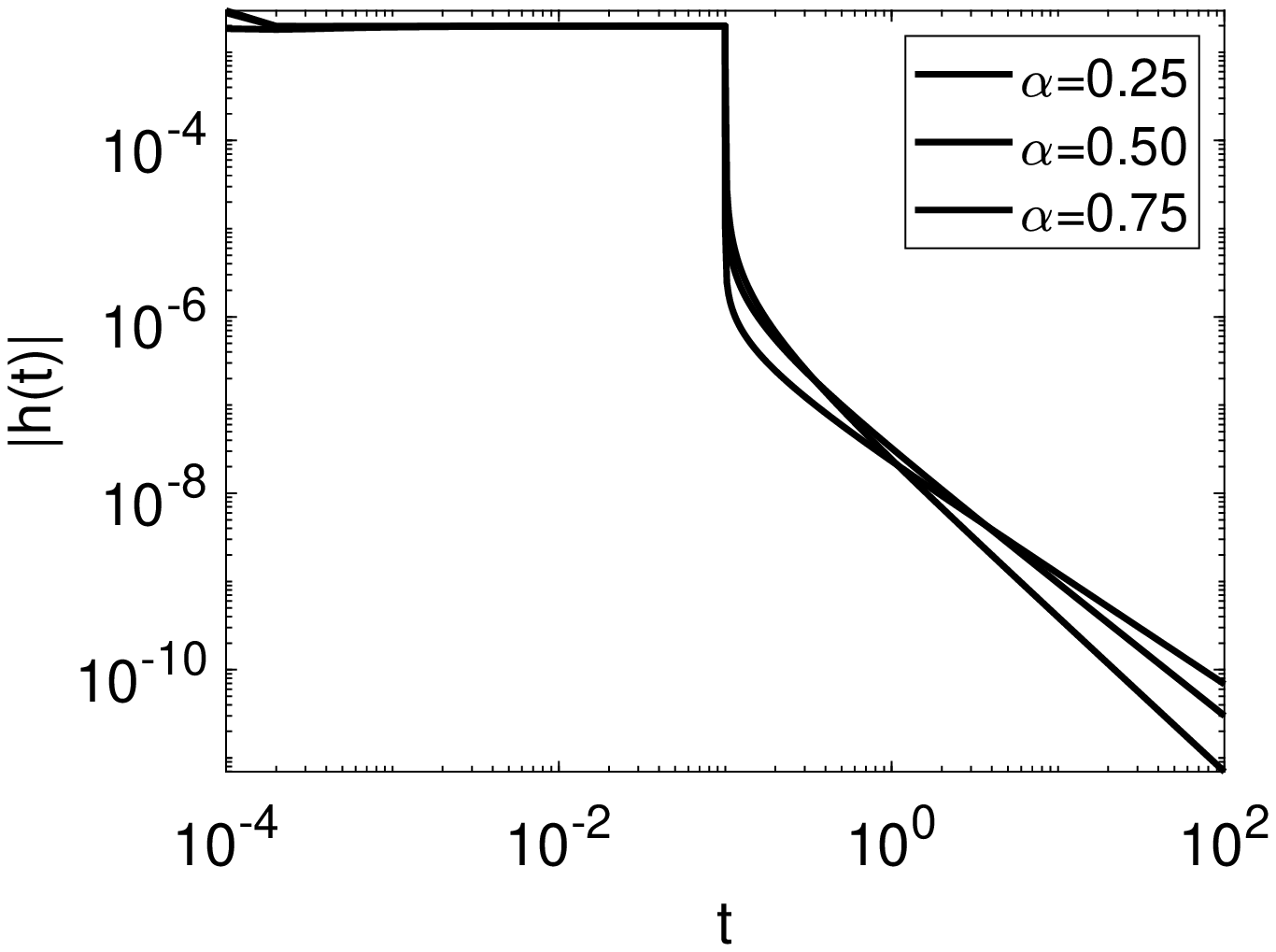} & \includegraphics[width=0.33\textwidth]{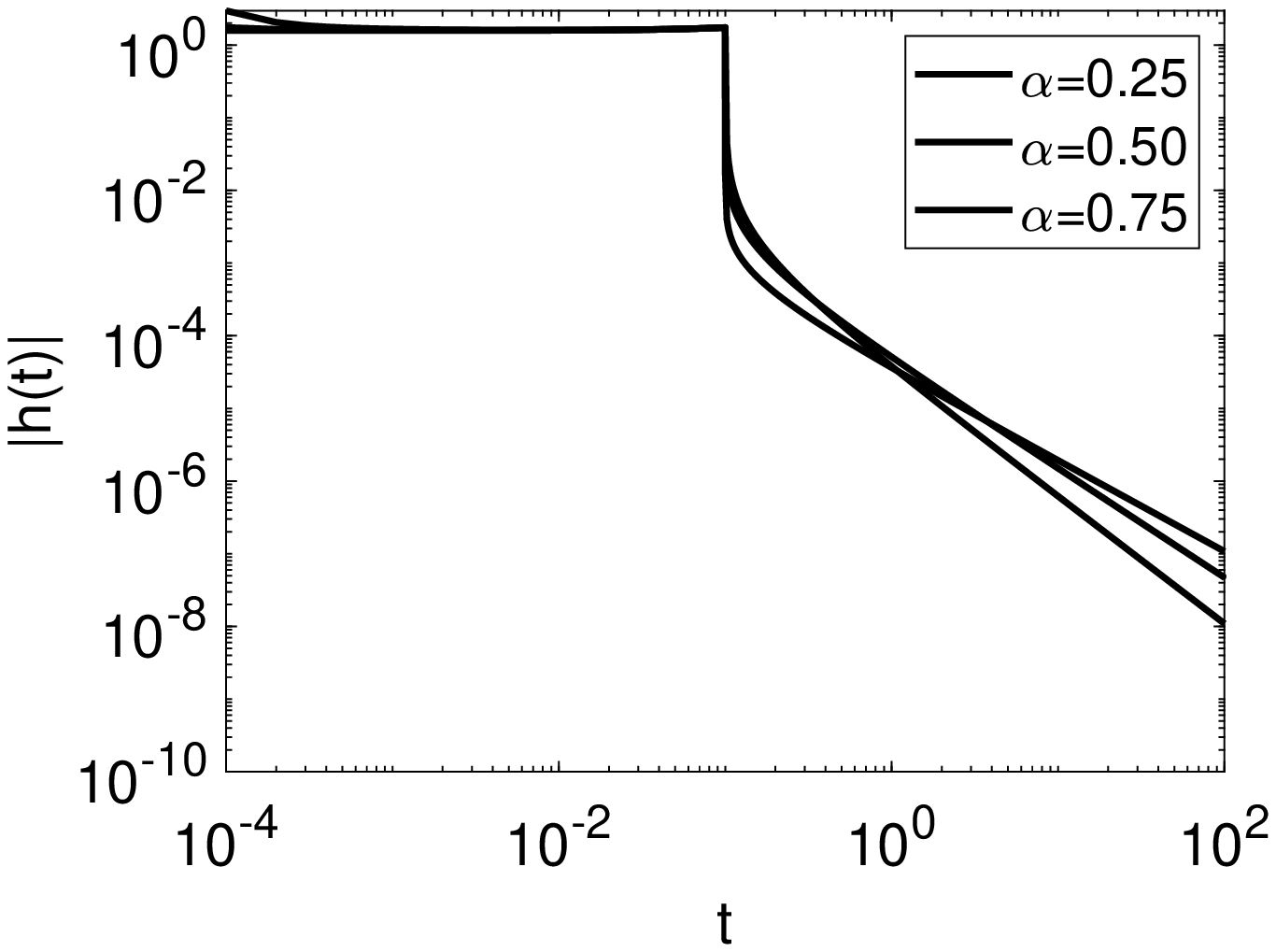}\\
   (i) &  (ii) &  (iii)
\end{tabular}
\caption{The profile of $|h(t)|$ for the three cases of Example \ref{exam:1dsub}.\label{fig:1dsub}}
\end{figure}

\begin{table}[htp!]
  \centering\small
  \begin{threeparttable}
  \caption{Recovery of the fractional order $\alpha$ for Example \ref{exam:1dsub}. The blocks (a), (b) and (c) are for the observation
  window $[1,2]$, $[1,10]$ and $[10,20]$, respectively.\label{tab:1dsub}}
    \begin{tabular}{ccccccccccc}
    \toprule
    \multicolumn{2}{c}{}&
    \multicolumn{3}{c}{(i)}&\multicolumn{3}{c}{(ii)} & \multicolumn{3}{c}{(iii)}\\
    \cmidrule(lr){3-5} \cmidrule(lr){6-8} \cmidrule(lr){9-11}
    & $\epsilon\backslash\alpha$ &0.25&0.5&0.75&0.25&0.50&0.75&0.25&0.50&0.75\\
    \midrule
       & 0\% & 0.249 & 0.500 & 0.750 & 0.297 & 0.557 & 0.817 & 0.298 & 0.558 & 0.818\\
   (a) & 1\% & 0.238 & 0.488 & 0.738 & 0.285 & 0.545 & 0.806 & 0.286 & 0.546 & 0.806\\
       & 5\% & 0.185 & 0.435 & 0.685 & 0.232 & 0.492 & 0.753 & 0.233 & 0.493 & 0.753\\
    \hline
       & 0\% & 0.249 & 0.500 & 0.750 & 0.273 & 0.528 & 0.783 & 0.273 & 0.528 & 0.783\\
   (b) & 1\% & 0.244 & 0.494 & 0.744 & 0.267 & 0.522 & 0.777 & 0.268 & 0.522 & 0.777\\
       & 5\% & 0.244 & 0.494 & 0.744 & 0.267 & 0.522 & 0.777 & 0.240 & 0.495 & 0.750\\
    \hline
       & 0\% & 0.249 & 0.500 & 0.750 & 0.254 & 0.505 & 0.756 & 0.254 & 0.505 & 0.756\\
   (c) & 1\% & 0.238 & 0.488 & 0.738 & 0.242 & 0.493 & 0.744 & 0.242 & 0.493 & 0.744\\
       & 5\% & 0.238 & 0.488 & 0.738 & 0.242 & 0.493 & 0.744 & 0.189 & 0.440 & 0.691\\
    \bottomrule
\end{tabular}
\end{threeparttable}
\end{table}

The profiles of the Neumann trace data ${h}(t)=\partial_\nu u(x_0,t)$ are shown in Fig. \ref{fig:1dsub}
(in the doubly logarithmic scale). Clearly, a power type decay is observed for large $t$ and the
decay is faster when the initial condition $u_0$ vanishes identically. This observation agrees well with
the theoretical analysis in Section \ref{sec:prelim}, cf. Propositions \ref{prop:asympt-u0}--\ref{prop:asympt-g}.
In particular, it indicates that by fitting fractional powers to the discrete observation
points, one may obtain reasonable estimate on the fractional order $\alpha$.

In Table \ref{tab:1dsub} we show the recovered order $\alpha$ for three different observation windows
$[T_1,T_2]$, i.e., $[1,2]$, $[1,10]$, and $[20,20]$. The results are obtained using one single term in the
least-squares formulation \eqref{eqn:lsq}. It is observed that both observation window $[T_1,T_2]$
and the accuracy of the data influence the quality of order recovery, and the
behavior is more or less just as expected: the recovered order $\alpha$ becomes less accurate as
the observation window size becomes smaller or the data $h^\delta$ gets noisier. When the window
$[T_1,T_2]$ is sufficiently large, the recovery procedure is stable, and can yield accurate results
for up to $5\%$ noise in the data. Generally,
the results for Case (i) are more accurate than that in Cases (ii) and (iii), indicating that
nonzero initial data excitation is preferred for order recovery. One surprising phenomenon in the
presence of data noise, the recovery accuracy can improve over exact data, when
only the source $F$ or the boundary data $h$ is nonvanishing. The mechanism of this observation remains unclear.
Moreover, as the theory predicts, the results for Cases (ii) and (iii) are close to each other.
These results show the feasibility of the recovering the order $\alpha$
without knowing the medium.

The next example is about two-dimensional subdiffusion on a square domain with a circular
inclusion, where $B_r({x})$ denotes a ball centered at $x$ with a radius $r$.
\begin{example}\label{exam:2dsub}
The domain $\widetilde \Omega=(0,1)^2$ and the observation point $x_0$ is $(0,0.5)$.
\begin{itemize}
\item[{\rm(i)}] $(\omega,a,q,u_0,F,g)=(B_{0.2}(0.5,0.5),1+\sin(\pi x_1)x_2(1-x_2),1, x_1(1-x_1)\sin(\pi x_2), x_1(1-x_1)x_2(1-x_2)t\chi_{[0,0.1]}(t),0)$.
\item[{\rm(ii)}] $(\omega,a,q,u_0,F,g)=(\emptyset,1,1,0,\sin(\pi x_1)x_2^2(1-x_2)\chi_{[0,0.1]}(t),0)$.
\item[{\rm(iii)}] $(\omega,a,q,u_0,F,g) = (B_{0.2}(0.5,0.5),1+\sin(\pi x_1)\sin(\pi x_2),1,0,0, x_1(1-x_1)e^t\chi_{[0,0.1]}(t))$, where the Dirichlet boundary condition $g$
is specified only on the bottom boundary $\{(x_1,0):0\leq x_1\le 1\}$, and zero else where.
\end{itemize}
\end{example}

\begin{figure}[hbt!]
  \centering
  \setlength{\tabcolsep}{0pt}
  \begin{tabular}{ccc}
    \includegraphics[width=0.33\textwidth]{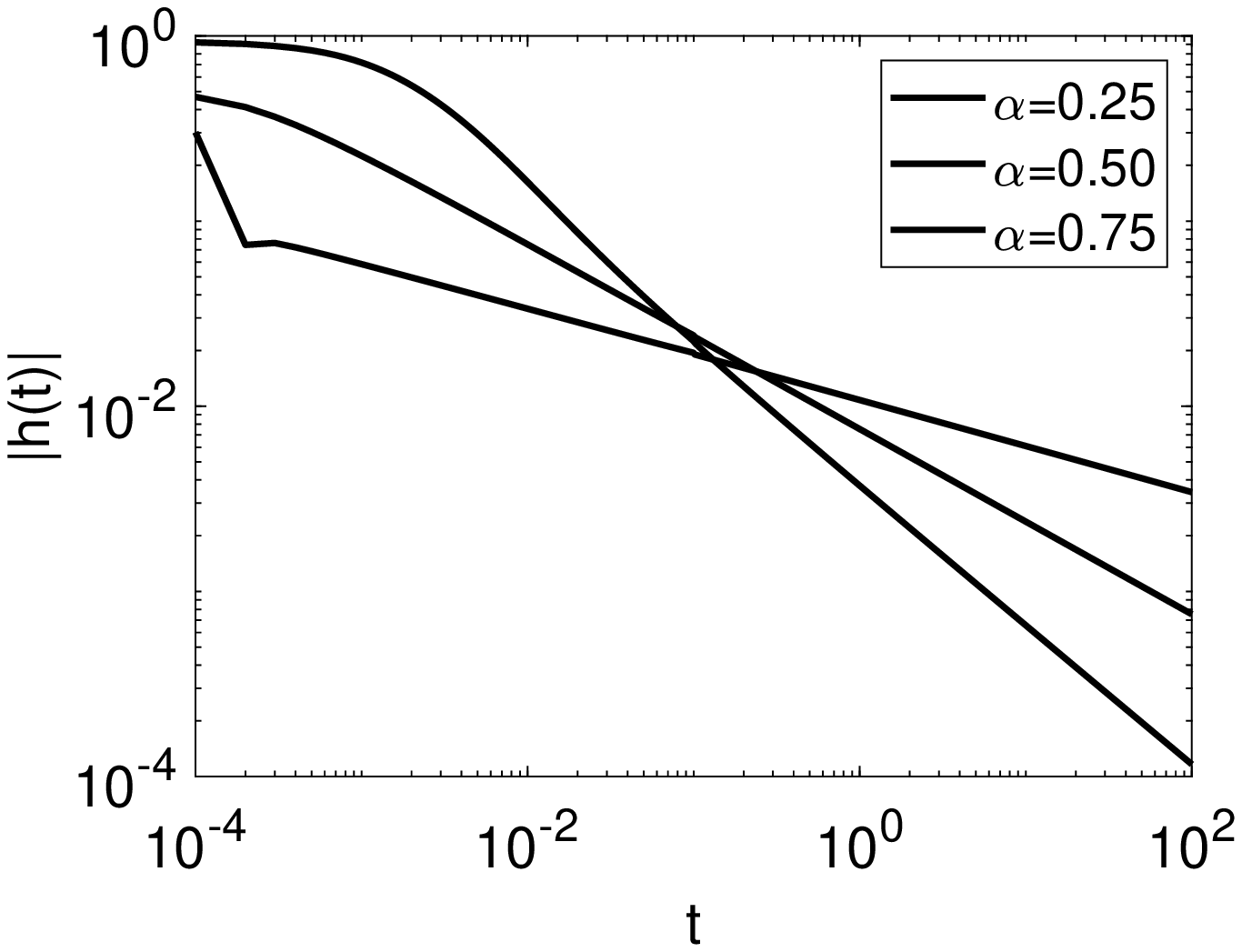} & \includegraphics[width=0.33\textwidth]{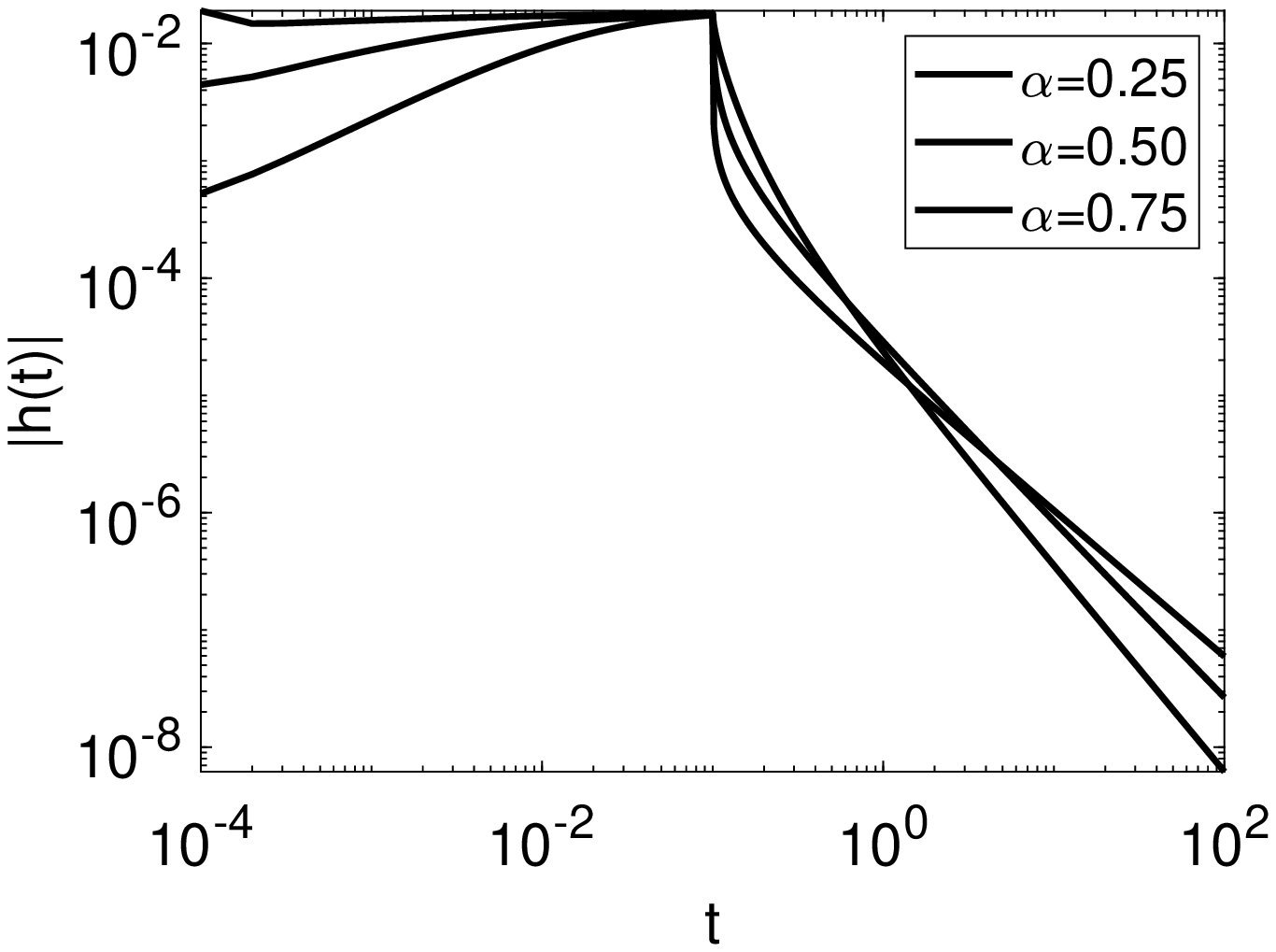} & \includegraphics[width=0.33\textwidth]{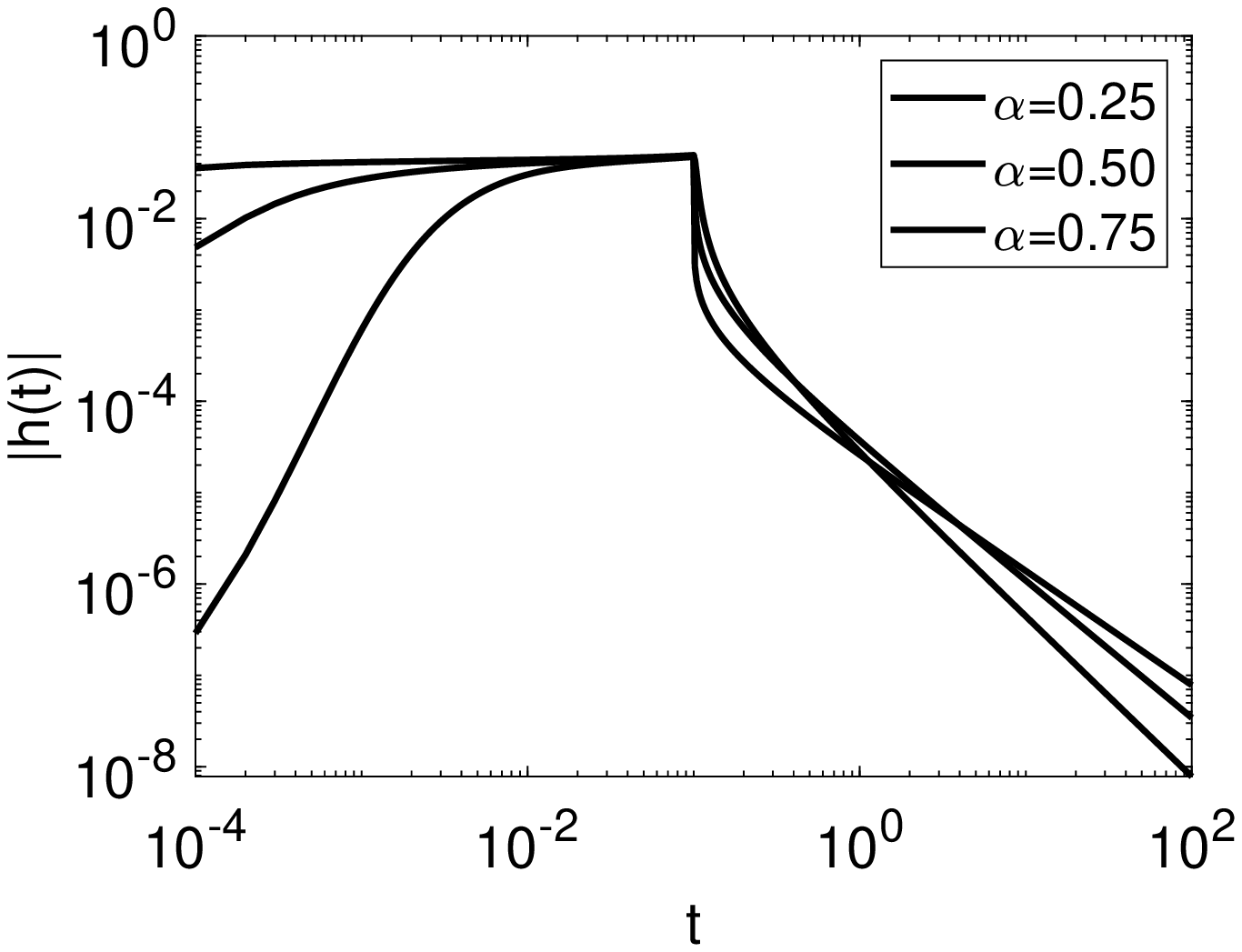}\\
    (i) & (ii) & (iii)
  \end{tabular}
  \caption{The profile of of $|{h}(t)|$ for the three cases of Example \ref{exam:2dsub}.\label{fig:2dsub}}
\end{figure}

The numerical results for Example \ref{exam:2dsub} are presented in Fig. \ref{fig:2dsub}
and Table \ref{tab:2dsub}. The decay behavior of the flux $h(t)=\partial_\nu u(x_0,t)$
is largely comparable with that for Example \ref{exam:1dsub}: after an initial transient
period, which is relatively short, the flux $h(t)$ shows a clearly power type decay, and
the decay is faster for cases (ii) and (iii) than case (i), confirming the theoretical
predictions from Propositions \ref{prop:asympt-u0}--\ref{prop:asympt-g}. The accuracy of the recovery is also
comparable with the one-dimensional case in Example \ref{exam:1dsub}. Note that the
presence of an obstacle $\omega$ within the domain $\Omega$ does not influence
much the recovery accuracy of the order $\alpha$, which agrees with the theoretical analysis.

\begin{table}[htp!]
  \centering\small
  \begin{threeparttable}
  \caption{Recovery of the fractional order $\alpha$ for Example \ref{exam:2dsub}. The blocks (a), (b) and (c) are for the observation
  windows $[1,2]$, $[1,10]$ and $[10,20]$, respectively.\label{tab:2dsub}}
    \begin{tabular}{ccccccccccc}
    \toprule
    \multicolumn{2}{c}{}&
    \multicolumn{3}{c}{(i)}&\multicolumn{3}{c}{(ii)} & \multicolumn{3}{c}{(ii)}\\
    \cmidrule(lr){3-5} \cmidrule(lr){6-8} \cmidrule(lr){9-11}
    & $\epsilon\backslash\alpha$ &0.25&0.5&0.75&0.25&0.50&0.75&0.25&0.50&0.75\\
    \midrule
       & 0\% & 0.248 & 0.500 & 0.756 & 0.281 & 0.554 & 0.879 & 0.293 & 0.558 & 0.837\\
   (a) & 1\% & 0.236 & 0.488 & 0.745 & 0.270 & 0.543 & 0.867 & 0.281 & 0.546 & 0.825\\
       & 5\% & 0.183 & 0.435 & 0.692 & 0.217 & 0.489 & 0.814 & 0.228 & 0.493 & 0.772\\
    \hline
       & 0\% & 0.248 & 0.500 & 0.754 & 0.260 & 0.526 & 0.818 & 0.269 & 0.528 & 0.794\\
   (b) & 1\% & 0.242 & 0.494 & 0.748 & 0.254 & 0.521 & 0.812 & 0.264 & 0.522 & 0.788\\
       & 5\% & 0.215 & 0.467 & 0.721 & 0.227 & 0.494 & 0.785 & 0.237 & 0.495 & 0.761\\
    \hline
       & 0\% & 0.249 & 0.500 & 0.751 & 0.245 & 0.505 & 0.766 & 0.251 & 0.505 & 0.759\\
   (c) & 1\% & 0.237 & 0.488 & 0.739 & 0.234 & 0.493 & 0.755 & 0.240 & 0.493 & 0.748\\
       & 5\% & 0.184 & 0.435 & 0.686 & 0.181 & 0.440 & 0.701 & 0.187 & 0.440 & 0.695\\
    \bottomrule
\end{tabular}
\end{threeparttable}
\end{table}
   \subsection{Numerical results for diffusion wave}
Now we present two-dimensional examples for the diffusion wave case. The setting is identical
with that of Example \ref{exam:2dsub} for subdiffusion, except the fractional order.
\begin{example}\label{exam:2ddw}
The domain $\widetilde \Omega=(0,1)^2$ and the observation point $x_0$ is $(0,0.5)$.
\begin{itemize}
\item[{\rm(i)}] $(\omega,a,q,u_0,F,g)=(B_{0.2}(0.5,0.5),1+\sin(\pi x_1)x_2(1-x_2),1, x_1(1-x_1)\sin(\pi x_2), x_1(1-x_1)x_2(1-x_2)t\chi_{[0,0.1]}(t),0)$.
\item[{\rm(ii)}] $(\omega,a,q,u_0,F,g)=(\emptyset,1,1,0,\sin(\pi x_1)x_2^2(1-x_2)\chi_{[0,0.1]}(t),0)$.
\item[{\rm(iii)}] $(\omega,a,q,u_0,F,g)=(B_{0.2}(0.5,0.5),1+\sin(\pi x_1)\sin(\pi x_2),1,0,0,x_1(1-x_1)e^t\chi_{[0,0.1]}(t))$, where the Dirichlet
boundary condition $g$ is specified only on the bottom boundary $\{(x_1,0): 0\leq x_1\le 1\}$, and zero else where.
\end{itemize}
\end{example}

\begin{figure}[hbt!]
  \centering
  \setlength{\tabcolsep}{0pt}
  \begin{tabular}{ccc}
    \includegraphics[width=0.33\textwidth]{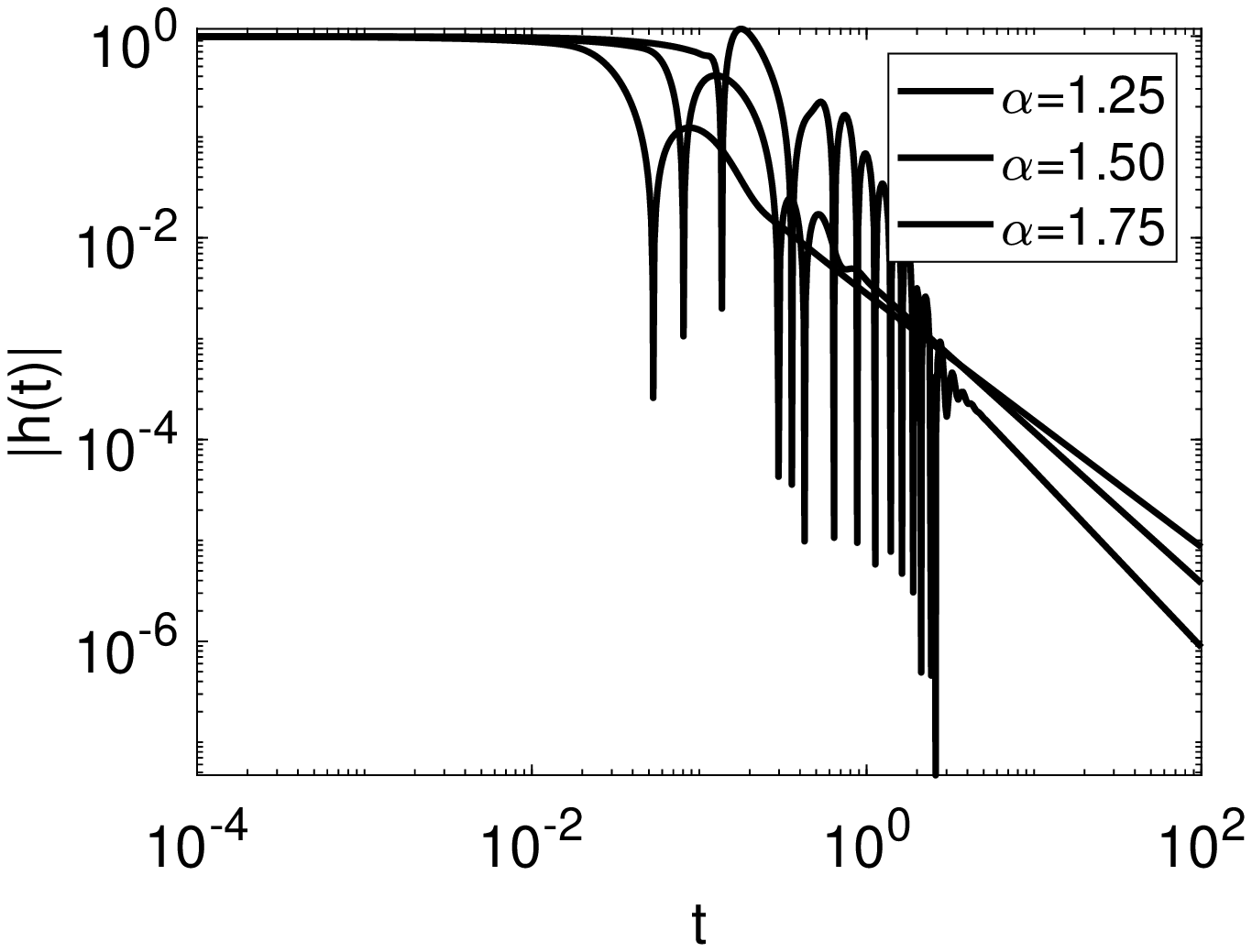} & \includegraphics[width=0.33\textwidth]{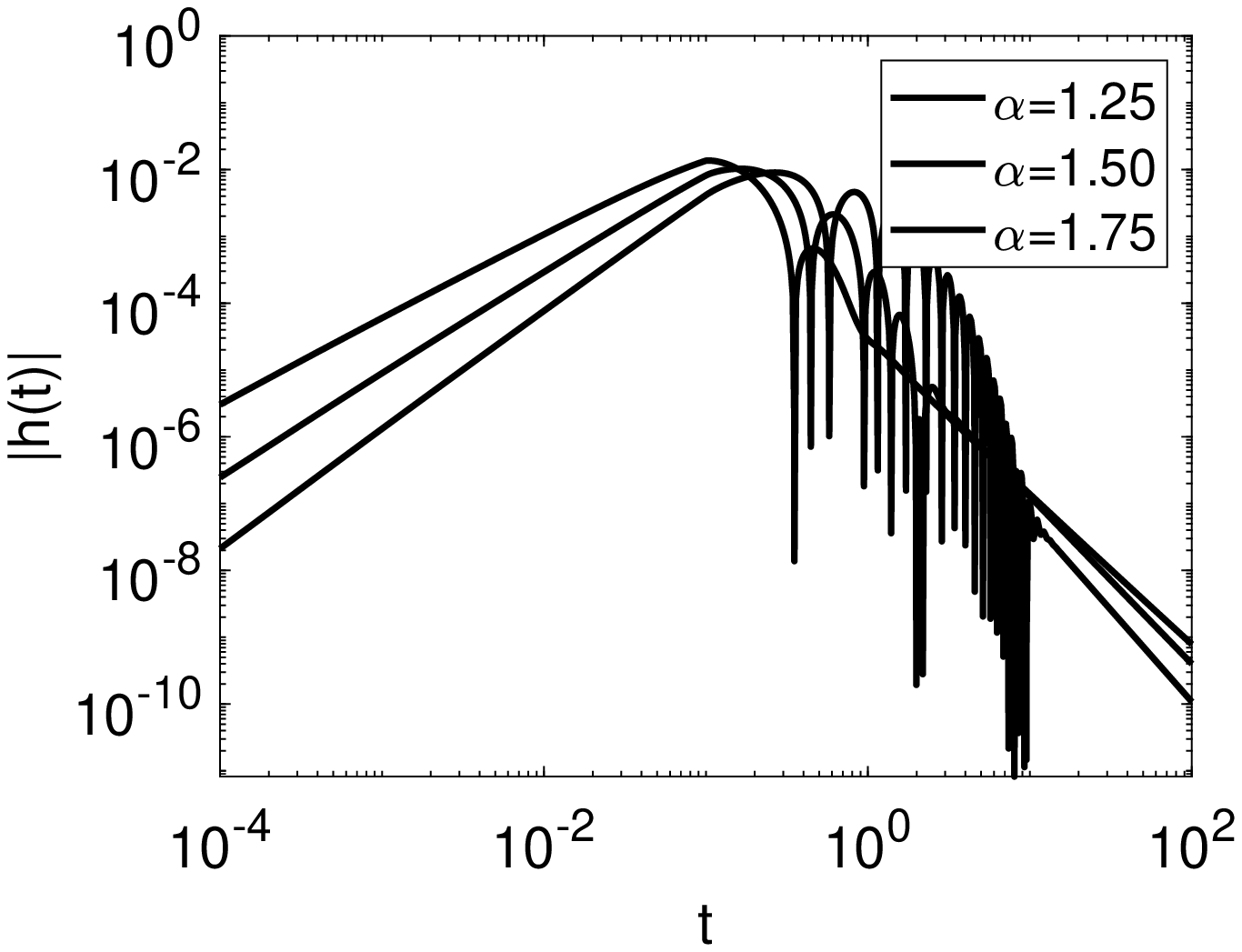} & \includegraphics[width=0.33\textwidth]{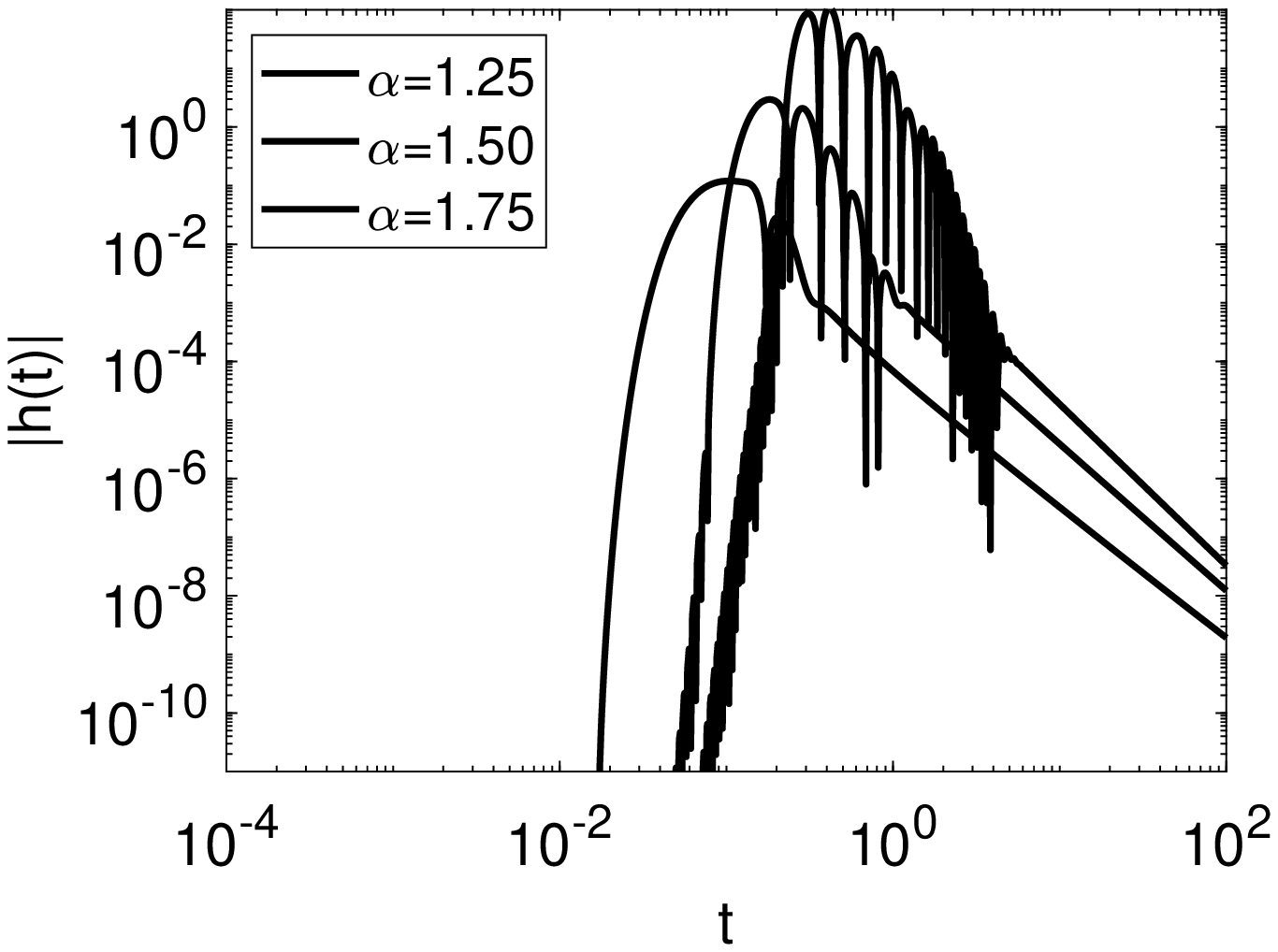}\\
    (i) & (ii) & (iii)
  \end{tabular}
  \caption{The profile of of $|h(t)|$ for the three cases of Example \ref{exam:2ddw}.\label{fig:2ddw}}
\end{figure}

\begin{table}[htp!]
  \centering
  \begin{threeparttable}
  \caption{Recovery of the fractional order $\alpha$ for Example \ref{exam:2ddw}. The blocks (a), (b) and (c) are for the observation
  windows $[1,10]$, $[20,30]$ and $[20,23]$, respectively.\label{tab:2ddw}}
    \begin{tabular}{ccccccccccc}
    \toprule
    \multicolumn{2}{c}{}&
    \multicolumn{3}{c}{(i)}&\multicolumn{3}{c}{(ii)} & \multicolumn{3}{c}{(iii)}\\
    \cmidrule(lr){3-5} \cmidrule(lr){6-8}\cmidrule(lr){9-11}
    & $\epsilon\backslash\alpha$ &1.25&1.50&1.75&1.25&1.50&1.75 &1.25&1.50&1.75\\
    \midrule
       & 0\% & 1.258 & 1.501 & --    & 1.321 & --    & --    & 1.315 & 1.567 & --   \\
   (a) & 1\% & 1.253 & 1.495 & --    & 1.315 & --    & --    & 1.309 & 1.561 & --   \\
       & 5\% & 1.226 & 1.468 & --    & 1.288 & --    & --    & 1.282 & 1.534 & --   \\
    \hline
       & 0\% & 1.250 & 1.500 & 1.749 & 1.251 & 1.504 & 1.752 & 1.238 & 1.503 & 1.754\\
   (b) & 1\% & 1.229 & 1.479 & 1.728 & 1.230 & 1.483 & 1.731 & 1.217 & 1.483 & 1.733\\
       & 5\% & 1.134 & 1.383 & 1.633 & 1.135 & 1.388 & 1.636 & 1.121 & 1.387 & 1.638\\
    \hline
       & 0\% & 1.250 & 1.500 & 1.749 & 1.253 & 1.505 & 1.752 & 1.241 & 1.504 & 1.755\\
   (c) & 1\% & 1.187 & 1.437 & 1.686 & 1.191 & 1.442 & 1.689 & 1.178 & 1.442 & 1.692\\
       & 5\% & 0.899 & 1.148 & 1.398 & 0.902 & 1.154 & 1.401 & 0.890 & 1.153 & 1.404\\
    \bottomrule
\end{tabular}
\end{threeparttable}
\end{table}
The profiles of the Neumann trace ${h}(t)=\partial_\nu u(x_0,t)$ are shown in Fig. \ref{fig:2ddw}.
Compared with the subdiffusion case, the trace $|h(t)|$ exhibits many more oscillations (or equivalently
$h(t)$ oscillates more widely around zero), and as a result, the transient period is much
longer. This behavior seems characteristic of the diffusion wave
problem: for $\alpha\in(1,2)$, the Mittag-Leffler functions $E_{\alpha,2}(-t)$ and $E_{\alpha,\alpha}
(-t)$ are no longer completely monotone, which is in stark contrast with that for the subdiffusion case (for
which both are completely monotone \cite{Pollard:1948,Schneider:1996} and thus does not change sign).
Further, the number of real roots of both functions increases to infinity as the order $\alpha$ tends to two;
see the work \cite{JinRundell:2012} for an empirical study on the roots of the function $E_{\alpha,2}
(-t)$. The plots in the middle and right panels show far more oscillations than that in the
left most panel (when the $\alpha$ value is the same). This might be related to the empirical observation
that for any fixed $\alpha\in (1,2)$, the function $E_{\alpha,\alpha}(-t)$ has more real roots than
$E_{\alpha,2}(-t)$ (which, however, has not been rigorously proved so far). Note that the magnitude of
$h(t)$ in Case (iii) is very small during the initial time, and thus not displayed in the plot, which
differs greatly from the subdiffusion case. In sum, in the diffusion wave case, the boundary data $h$ does exhibits
a power type decay for large time $t$, but the asymptotic
power decay kicks in only for much larger $t$, which is especially pronounced for the order $\alpha$
close to two. These observations necessitate measurements at large time so that the least-squares formulation
\eqref{eqn:lsq} is numerically viable.

The numerical recovery results for Example \ref{exam:2ddw} are given in Table \ref{tab:2ddw}. Just as the plots in
Fig. \ref{fig:2ddw} predict, when the initial time $T_1$ of the observation window $[T_1,T_2]$ is not sufficient
large, the least-squares approach fails to deliver reasonable recovery for large $\alpha$, as is indicated by notation
``--'' in the table. This is more dramatic for Cases (ii) and (iii) than Case (i), and it is attributed to the fact
when $T_1$ is small, the data is still too far away from the asymptotic regime on which the least-squares formulation
\eqref{eqn:lsq} is based. When the initial time $T_1$ of the window $[T_1,T_2]$ increases, the recovery becomes viable
again and the recovered orders are accurate for data with up to $5\%$ noise, indicating the necessity of large
initial observation time $T_1$. When the window size decreases from ten to three, the stability of the recovery
worsens quite a bit, as confirmed by the numerical results in blocks (b) and (c) in Table \ref{tab:2ddw}.

\appendix
\section{Numerical schemes for the direct problem \eqref{eq1}}
In this appendix, we describe the numerical schemes for simulating the direct problem \eqref{eq1} for completeness;
see the review \cite{JinLazarovZhou:2019cmame} for further details. For the spatial discretization,
we employ the Galerkin finite element method with continuous piecewise linear finite element basis.
Let $X_h$ be the continuous piecewise linear finite element space, subordinated to a shape regular
triangulation of the domain $\Omega$, and $M_h$ and $S_h$ be the corresponding mass and stiffness
matrices, respectively. The temporal discretization is based on the finite-difference approximation.
For any $N\in\mathbb{N}$ total number of time steps, let
$\tau=\frac{T}{N}$ be the time step size, and $t_n=n\tau$, $n=0,\ldots, N$, the time grid. We define
the difference approximations (with the shorthand $u^n=u(t_n)$)
$\delta_tu^{n+\frac12}=\tau^{-1}({u^{n+1}-u^n})$ and $\delta_t^2u^n=\tau^{-1}(\delta_tu^{n+\frac12}-\delta_tu^{n-\frac12})$.

Note that a direct implementation of many time stepping schemes suffers from a serious storage issue, due to
the nonlocality of the operator $\partial_t^\alpha u$. Below we describe an implementation based on the sum of exponentials
(SOE) approximation of the function $t^{-\beta}$ over a compact interval $[\delta,T]$ (with
$\delta>0$) \cite{Beylkin:2005,Jiang:2017,McLean:2018}. In practice, with proper model reduction,
tens of terms suffice a reasonable accuracy.
\begin{lemma}\label{prop:soe}
For any $0<\beta<2$, $0<\delta<1$ and $0<\varepsilon<1$, there exist $\{(s_i,\omega_i)\}_{i=1}^{N_e}\subset\mathbb{R}_+^2$ such that
$\Big|t^{-\beta}-\sum_{i=1}^{N_e}\omega_ie^{-s_it}\Big|\leq \varepsilon$, for all $t\in [\delta,T],$
with $N_{e}=\mathcal{O}((\log \frac1\varepsilon)(\log\log\frac1\varepsilon+\log\frac{T}{\delta})+(\log\frac1\delta)(\log\log\frac1\varepsilon+\log\frac1\delta))$.
\end{lemma}

\subsection{Numerical scheme for subdiffusion ($0<\alpha<1$)}
Using piecewise linear interpolant, integration by parts, and the SOE approximation (with $\beta=1+\alpha$
and $\delta=\tau$), we can approximate the Djrbashian-Caputo fractional derivative $\partial_t^\alpha u(t_n)$, $n\geq1$, by
\begin{align*}
  \partial_t^\alpha u^n &= \frac{1}{\Gamma(1-\alpha)}\int_{t_{n-1}}^{t_n}(t_n-s)^{-\alpha}u'(s)\d s +\frac{1}{\Gamma(1-\alpha)}\int_0^{t_{n-1}}(t_n-s)^{-\alpha}u'(s)\d s\\
    & \approx \frac{u^n-u^{n-1}}{\tau^\alpha c_\alpha} + \frac{1}{\Gamma(1-\alpha)}\Big[\frac{u^{n-1}}{\tau^\alpha}-\frac{u^0}{t_n^\alpha} - \alpha\sum_{i=1}^{N_{e}}\omega_iU_{h,i}^n\Big],
\end{align*}
with the history terms
$$U_{h,i}^n = \int_0^{t_{n-1}}e^{-(t_n-s)s_i}u(s)\d s,$$
and $c_\alpha=\Gamma(2-\alpha)$. Since $u(s)$ is piecewise
linear (i.e., $u(s)=u^{n-2}+\delta_tu^{n-\frac32}(s-t_{n-2})$ over $[t_{n-2},t_{n-1}]$), $U_{h,i}^n$ satisfies the following recursion
\begin{align*}
  U_{h,i}^n &= e^{-s_i\tau} U_{h,i}^{n-1} + \int_{t_{n-2}}^{t_{n-1}}e^{-s_i(t_n-s)}u(s)\d s
       = e^{-s_i\tau} U_{h,i}^{n-1} + w_i^1u^{n-1}+w_i^2u^{n-2},
\end{align*}
with the weights
$$w_i^1=\frac{e^{-s_i\tau}}{s_i^2\tau}(e^{-s_i\tau}-1+s_i\tau)\quad\mbox{and}\quad
w_i^2=\frac{e^{-s_i\tau}}{s_i^2\tau}(1-e^{-s_i\tau}-e^{-s_i\tau}s_i\tau).$$
When $s_i\tau$ is very small, the computation of the weights $w_i^1$ and $w_i^2$ is prone to
cancellation errors. Then they may be computed by Taylor expansion as
$$w_i^1 \approx e^{-s_i\tau}\tau\Big(\frac{1}{2}-\frac{s_i\tau}{6}+\frac{(s_i\tau)^2}{24}\Big)\quad\mbox{and}\quad
w_i^2 \approx e^{-s_i\tau}\tau\Big(\frac{1}{2}-\frac{s_i\tau}{3} + \frac{(s_i\tau)^2}{8}\Big).$$
The fully discrete scheme reads: with $U_h^0=P_hu_0$ (with $P_h$ being the $L^2$ projection on
$X_h$), find $U_h^i\in X_h$ such that for $n=1,2,\ldots,N$,
\begin{align*}
   (M_h+c_\alpha\tau^\alpha S_h)U_h^n = \alpha M_h U_h^{n-1} + (1-\alpha)M_h\Big[n^{-\alpha}U_h^0+\alpha\tau^\alpha\sum_{i=1}^{N_e}\omega_iU_{h,i}^n\Big] + c_\alpha\tau^\alpha F_h^n.
\end{align*}
This scheme has an accuracy $O(\tau^{2-\alpha})$ for smooth solutions \cite{LinXu:2007} and
$O(\tau)$ for general incompatible problem data \cite{JinLazarovZhou:2016ima}. The first
step may be corrected to be  \cite{YanKhanFord:2018}
\begin{equation*}
  (M_h+c_\alpha\tau^\alpha S_h)U_h^1 = M_h U_h^0 + c_\alpha\tau^\alpha(F_h^1 + \tfrac{1}{2}F_h^0-\tfrac{1}{2}S_hU_h^0).
\end{equation*}
Then the overall accuracy of the corrected scheme is $O(\tau^{2-\alpha})$ \cite{YanKhanFord:2018}.

\subsection{Numerical scheme for diffusion wave ($1<\alpha<2$)}
For $1<\alpha<2$, we employ piecewise quadratic interpolant of $u$. Let $H_{2,0}(t)$
be the Hermit quadratic interpolant through $(0,u^0), (\tau,u^1)$ and $(0,u^{\prime 0})$:
$$H_{2,0}(t)=u^0+u^{\prime0}t + \tau^{-1}(\delta_t u^\frac12-u^{\prime0})t^2,$$
with $H_{2,0}''(t)=2\tau^{-1}(\delta_tu^\frac12-u'(t_0))$. Next, for any $u$ defined on the
interval $[t_{n-1},t_{n+1}]$, $n=1,\ldots, N-1$, using $(t_{n-1},u^{n-1})$, $(t_n,u^n)$, $(t_{n+1},u^{n+1})$, let
$L_{2,n}(t)$ be the quadratic Lagrangian interpolant
$$L_{2,n}(t)=u^{n-1} + (\delta_tu^{n-\frac12})(t-t_{n-1}) + \tfrac12(\delta_t^2u^n)(t-t_{n-1})(t-t_n),$$ with
 $L_{2,n}''(t)=\delta_t^2u^n$. The scheme employs $H_{2,0}(t)$ on $[t_0,t_{\frac12}]$ and $L_{2,n}$
on $[t_{n-\frac12},t_{n+\frac12}]$, and SOE approximation
(with $\beta=\alpha-1$ and $\delta=\frac{\tau}{2}$). Then for $n=1$, there holds
\begin{align*}
   \partial_t^\alpha u^\frac12 &= \frac{1}{c_\alpha}\int_0^{t_\frac12}(t_\frac12-s)^{1-\alpha}u''(s)\d s
   \approx \frac{1}{c_\alpha}\int_0^{t_\frac12}H_{2,0}''(s)(t_\frac12-s)^{1-\alpha}\d s
      =\frac{2^{\alpha-1}(u^1-u^0-\tau u^{\prime0})}{\tau^{\alpha}c_\alpha'},
\end{align*}
with $c_\alpha=\Gamma(2-\alpha)$ and $c_\alpha'=\Gamma(3-\alpha)$. Similarly, for $n=2$, we have
\begin{align*}
   \partial_t^\alpha u^\frac32&= \frac{1}{c_\alpha}\Big[\int_0^{t_\frac12}(t_{n-\frac12}-s)^{1-\alpha}u''(s)\d s+\int_{t_{\frac12}}^{t_{\frac32}}(t_{\frac32}-s)^{1-\alpha}u''(s)\d s\Big]\\
   &\approx \frac{1}{c_\alpha}\Big[\int_0^{t_\frac12}H_{2,0}''(s)\sum_{i=1}^{N_e}\omega_ie^{-s_i(t_{n-\frac12}-s)}\d s + \int_{t_{\frac12}}^{t_{\frac32}}L_{2,1}''(s)(t_{\frac32}-s)^{1-\alpha}\d s \Big]\\
   & = \frac{1}{c_\alpha}\sum_{i=1}^{N_e}\omega_iU_{h,i}^2 + \frac{\delta_t^2u^1}{c_\alpha'}\tau^{2-\alpha}= \frac{1}{c_\alpha}\sum_{i=1}^{N_e}\omega_iU_{h,i}^2 + \frac{u^2-2u^1+u_0}{c_\alpha'\tau^\alpha},
\end{align*}
with the history term $U_{h,i}^2$ given by
\begin{equation*}
  U_{h,i}^2 = \int_{0}^{t_\frac12}H_{2,0}''(s)e^{-s_i(t_\frac32-s)}\d s = \frac{2(\delta_t u^{\frac12}-u^{\prime0})}{s_i\tau}e^{-s_i\tau}(1-e^{-\frac{1}{2}s_i\tau}),\quad i=1,\ldots,N_e.
\end{equation*}
For small $s_i\tau$,
$$\frac{2(1-e^{-\frac12s_i\tau})}{s_i\tau}\approx1-\frac{s_i\tau}{4}+\frac{(s_i\tau)^2}{24}.$$
Last, for $n\geq3$,
\begin{align*}
   &\quad\partial_t^\alpha u(t_{n-\frac12})= \frac{1}{c_\alpha}\Big[\int_0^{t_\frac12}(t_{n-\frac12}-s)^{1-\alpha}u''(s)\d s+\sum_{k=1}^{n-1}\int_{t_{k-\frac{1}{2}}}^{t_{k+\frac12}}(t_{n-\frac12}-s)^{1-\alpha}u''(s)\d s\Big]\\
   &\approx \frac{1}{c_\alpha}\Big[\int_0^{t_\frac12}H_{2,0}''(s)\sum_{i=1}^{N_e}\omega_ie^{-s_i(t_{n-\frac12}-s)}\d s
   + \sum_{k=1}^{n-2}\int_{t_{k-\frac{1}{2}}}^{t_{k+\frac12}}L_{2,k}''(s)\sum_{i=1}^{N_e}\omega_ie^{-s_i(t_{n-\frac12}-s)}\d s\\
   &\quad + \int_{t_{n-\frac{3}{2}}}^{t_{n-\frac12}}L_{2,n-1}''(s)(t_{n-\frac12}-s)^{1-\alpha}\d s \Big]\\
   & = \frac{1}{c_\alpha}\sum_{i=1}^{N_e}\omega_iU_{h,i}^n + \frac{\delta_t^2u^{n-1}}{c_\alpha'}\tau^{2-\alpha}
   = \frac{1}{c_\alpha}\sum_{i=1}^{N_e}\omega_iU_{h,i}^n + \frac{u^n-2u^{n-1}+u^{n-2}}{\tau^{\alpha}c_\alpha'},
\end{align*}
where the history term $U_{h,i}^n$ is given by
\begin{equation*}
  U_{h,i}^n = \int_0^{t_\frac12}H_{2,0}''(s)e^{-s_i(t_{n-\frac12}-s)}\d s
   + \sum_{k=1}^{n-2}\int_{t_{k-\frac{1}{2}}}^{t_{k+\frac12}}L_{2,k}''(s)e^{-s_i(t_{n-\frac12}-s)}\d s.
\end{equation*}
Note that the history term $U_{h,i}^n$, $n=3,4,\ldots,$ can be evaluated recursively as
\begin{align*}
  U_{h,i}^n 
  & = e^{-s_i\tau}U_{h,i}^{n-1} + \delta_t^2 u^{n-2} \int_{t_{n-\frac{5}{2}}}^{t_{n-\frac32}}e^{-s_i(t_{n-\frac12}-s)}\d s
   = e^{-s_i\tau}U_{h,i}^{n-1} +\delta_t^2u^{n-2}\frac{e^{-s_i\tau}}{s_i}(1-e^{-s_i\tau}).
\end{align*}
For small $s_i\tau$, the weight $$\frac{1-e^{-s_i\tau}}{s_i\tau}\approx1-\frac{s_i\tau}{2}+\frac{(s_i\tau)^2}{6}-\frac{(s_i\tau)^3}{24}.$$
Then applying the Galerkin finite element method in space, we obtain
\begin{align*}
  M_hU^1 + \tfrac{c_\alpha'}{2^{\alpha}}\tau^\alpha S_h U^1 &= M_h(U^0+\tau U^{\prime0})-\tfrac{c_\alpha'}{2^{\alpha}}\tau^\alpha S_hU^0 + \tfrac{c_\alpha'}{2^{\alpha-1}}\tau^\alpha F_h^\frac{1}{2},\\
  M_hU^n + \tfrac{c_\alpha'}{2}\tau^\alpha S_h U^n &= M_h(2U^{n-1}-U^{n-2}) - \tfrac{c_\alpha'}{2}\tau^\alpha S_hU^{n-1} - (2-\alpha)\tau^\alpha M_h\sum_{i=1}^{N_e}\omega_iU_{h,i}^n\\
  &\quad  + c_\alpha'\tau^\alpha F_h^{n-\frac12},\quad n=2,3,\ldots.
\end{align*}
This scheme is expected to be $O(\tau)$ accurate for general problem data.

\bibliographystyle{abbrv}
\bibliography{frac}

%
\end{document}